\documentclass[11pt]{article}

\usepackage[a4paper,margin=30mm]{geometry}
\usepackage[T1]{fontenc}
\usepackage[utf8]{inputenc}
\usepackage{lmodern}
\usepackage{amsmath,amssymb,amsthm}
\usepackage{needspace}
\usepackage{xcolor}
\usepackage[colorlinks=true,linkcolor=blue!55!black,citecolor=blue!55!black,urlcolor=blue!55!black]{hyperref}

\newtheorem{theorem}{Theorem}[section]
\newtheorem{proposition}[theorem]{Proposition}
\newtheorem{lemma}[theorem]{Lemma}
\newtheorem{corollary}[theorem]{Corollary}

\theoremstyle{definition}
\newtheorem{definition}[theorem]{Definition}
\newtheorem{remark}[theorem]{Remark}
\newtheorem{example}[theorem]{Example}

\newcommand{\C}{\mathbb{C}}
\newcommand{\R}{\mathbb{R}}
\newcommand{\B}{\mathbb{B}}
\newcommand{\Pp}{\mathbb{P}}
\newcommand{\cO}{\mathcal{O}}
\newcommand{\tX}{\widetilde X}
\newcommand{\sqrti}{\sqrt{-1}}
\newcommand{\ddbar}{\sqrti\,\partial\overline\partial}
\newcommand{\rank}{\operatorname{rank}}
\newcommand{\Aut}{\operatorname{Aut}}
\newcommand{\PU}{\mathrm{PU}}
\title{Stein Criteria for Holomorphic Fibre Bundles}
\author{Masakazu Takakura}
\date{}

\begin{document}

\maketitle

\begin{abstract}
In this paper, we study Steinness of holomorphic fibre bundles. Over a compact Riemann surface, we prove that every holomorphic ball bundle is weakly $1$-complete and that its total space is Stein if and only if it admits no holomorphic section. We also obtain higher-dimensional Stein criteria for reductive ball bundles and affine bundles in terms of positivity of the canonical class. A topological obstruction to the $1$-convexity of holomorphic fibre bundles is also established.
\end{abstract}

\setcounter{tocdepth}{2}
\tableofcontents

\section{Introduction and Main Results}\label{sec:introduction}

\subsection{Background}\label{subsec:background}

Let $X$ be a complex manifold. We say that $X$ is weakly $1$-complete if there exists a smooth exhaustion function $\varphi$ such that $\ddbar\,\varphi\geq 0$. We say that $X$ is $1$-complete if there exists a smooth exhaustion function $\varphi$ such that $\ddbar\,\varphi>0$. Grauert's solution of the Levi problem implies that a complex manifold admitting a smooth strictly plurisubharmonic exhaustion is Stein; see \cite[Theorem~2 and discussion on pp.~468--469]{Grauert1958}. Let $\Omega\subset X$ be a domain. We say that $\Omega$ is locally pseudoconvex if, for every $p\in\partial\Omega$, there exists an open neighbourhood $U$ of $p$ such that $U\cap\Omega$ is Stein.

These notions of pseudoconvexity have been studied since the early development of several complex variables. Oka solved the Levi problem for domains in $\C^m$; see \cite[pp.~97--155]{Oka1953}. In particular, for a domain in $\C^m$, local pseudoconvexity, weak $1$-completeness, and $1$-completeness are equivalent. Docquier and Grauert extended the Levi problem to domains over Stein manifolds; see \cite[pp.~94--123]{DocquierGrauert1960}.

A natural question is how large the gap between local pseudoconvexity and weak $1$-completeness can be on a general complex manifold. Diederich and Forn{\ae}ss gave the following example.

\begin{theorem}[Diederich--Forn{\ae}ss {\normalfont\cite[pp.~119--123]{DiederichFornaess1982}}]\label{thm:diederich-fornaess-intro}
Let $S$ be a Hopf surface. Then there exist a holomorphic $\B^1$-bundle $E\to S$ and a holomorphic $\Pp^1$-bundle $F\to S$ such that $E\subset F$ is a locally pseudoconvex domain with real-analytic boundary, but $E$ is not weakly $1$-complete.
\end{theorem}

Interestingly, the situation changes substantially when the base is K\"ahler. Diederich and Ohsawa proved the following.

\begin{theorem}[Diederich--Ohsawa {\normalfont\cite[Theorem~1]{DiederichOhsawa1985}}]\label{thm:diederich-ohsawa-intro}
Every locally trivial holomorphic $\B^1$-bundle over a compact K\"ahler manifold is weakly $1$-complete.
\end{theorem}

Under the additional assumption of reductive monodromy, Seo proved weak $1$-completeness for bundles with bounded symmetric-domain fibres; see \cite[\S4]{Seo2022}. We ask when such total spaces are Stein and seek criteria formulated in terms of K\"ahler geometry. The present paper treats ball bundles and affine bundles.

\subsection{Steinness of Ball Bundles over Compact K\"ahler Manifolds}\label{subsec:intro-ball}

Let $X$ be a compact K\"ahler manifold, let $\widetilde X\to X$ be its universal covering, and let $\rho\colon\pi_1(X)\to\Aut(\B^m)$ be a representation. We denote by $E_\rho:=\widetilde X\times_\rho\B^m$ the associated holomorphic ball bundle, where $(\widetilde x,z)\sim(\gamma\widetilde x,\rho(\gamma)z)$ for $\gamma\in\pi_1(X)$. For holomorphic $\B^1$-bundles over compact Riemann surfaces, Diederich--Ohsawa and Barrett proved that the total space is Stein exactly when it contains no positive-dimensional compact analytic subset, equivalently when it admits no holomorphic section; see \cite[Theorem~1]{DiederichOhsawa1985} and \cite[Theorem~2]{Barrett1992}. We obtain the following extension to ball bundles with higher-dimensional fibres. The proof uses Adachi--Lee--Seo's argument establishing $1$-convexity of ball bundles \cite[Theorem~1.2]{AdachiLeeSeo2026}.

\begin{theorem}\label{thm:riemann-ball-criterion}
Let $R$ be a connected compact Riemann surface and let $p\colon E\to R$ be a holomorphic fibre bundle with fibre $\B^m$, $m\ge1$. Then $E$ is weakly $1$-complete. Moreover, the following conditions are equivalent:
\begin{enumerate}
    \item $E$ is Stein;
    \item $p\colon E\to R$ admits no holomorphic section;
    \item $E$ contains no positive-dimensional compact analytic subset.
\end{enumerate}
\end{theorem}

However, the compact-analytic-subset criterion in Theorem~\ref{thm:riemann-ball-criterion} is special to one-dimen\-sional bases. Example~\ref{ex:torus-B1-bundle} gives a reductive holomorphic $\B^1$-bundle over a compact K\"ahler surface which contains no positive-dimensional compact analytic subset but is not Stein.

We now consider bases of complex dimension at least two. In the reductive case, the existence theorem of Corlette and the rigidity theorems of Siu, Sampson, and Carlson--Toledo are particularly useful; see \cite[Corollary~3.5]{Corlette1988}, \cite[Proposition~3, pp.~81--84]{Siu1980}, and \cite[Theorems~2.3 and~3.5, Corollary~3.7, pp.~180--183]{CarlsonToledo1989}. They provide a distinguished pluriharmonic section and force it into a small number of possibilities: it is holomorphic, anti-holomorphic, or has low rank, in which case it is controlled by lower-dimensional geometry (in real rank two, it locally factors through a Riemann surface). Thus, although the base may be geometrically complicated, the harmonic section is highly constrained. This leads to a Stein criterion specific to the higher-dimensional setting. Let $c_\rho$ denote the Toledo class introduced in Definition~\ref{def:toledo-class}.

\begin{theorem}\label{thm:main-ball}
Let $X$ be a compact K\"ahler manifold with $\dim_{\C}X\geq 2$, and let $\rho\colon\pi_1(X)\to\Aut(\B^m)$ be reductive. Then the following are equivalent:
\begin{enumerate}
    \item $E_\rho$ is Stein;
    \item $c_\rho$ is a K\"ahler class.
\end{enumerate}
\end{theorem}

Theorem~\ref{thm:main-ball} does not cover the case of Riemann surfaces. In that case there are Stein ball bundles with reductive monodromy and $c_\rho=0$; see Example~\ref{ex:curve-rankone-zero-toledo}, where the harmonic section has real rank one and its image lies in a real geodesic. Thus the Stein criteria differ between one-dimensional and higher-dimensional bases.

This criterion is analogous to the classical line-bundle picture. By the Kodaira embedding theorem \cite[Main Theorem, pp.~28--48]{Kodaira1954}, a line bundle $L$ is ample if and only if $c_1(L)$ is a K\"ahler class. By Grauert's criterion \cite[Satz~8]{Grauert1962}, the total space of the dual line bundle $L^*\to X$ is $1$-convex with exceptional zero section if and only if $L$ is ample. Thus, in both settings, analytic convexity of the total space is governed by positivity of a natural cohomology class on the base.

Adachi--Lee--Seo further conjectured that a $1$-convex holomorphic $\B^m$-bundle over a compact K\"ahler manifold $X$ must satisfy $m\ge\dim_{\C}X$ \cite[Conjecture~1.3]{AdachiLeeSeo2026}. The same dimension obstruction holds more generally for bundles with bounded symmetric-domain fibres.

\begin{theorem}\label{thm:bsd-oneconvex-dimension}
Let $X$ be a (not necessarily K\"ahler) compact complex manifold of dimension $n$, let $\Omega$ be a bounded symmetric domain of dimension $m$, and let $p\colon E\to X$ be a holomorphic fibre bundle with fibre $\Omega$. If $E$ is $1$-convex, then $m\ge n$.
\end{theorem}

For ball bundles, imposing in addition that the base be compact K\"ahler and that the holonomy representation be reductive yields a more precise description.

\begin{theorem}\label{thm:main-ball-oneconvex}
Let $X$ be a compact K\"ahler manifold with $\dim_{\C}X\geq 2$, and let $\rho\colon\pi_1(X)\to\Aut(\B^m)$ be reductive. Suppose that $E_\rho$ is $1$-convex. Then $E_\rho$ admits a section $s\colon X\to E_\rho$ which is either holomorphic or anti-holomorphic. Moreover, one of the following holds:
\begin{enumerate}
    \item If $s$ is anti-holomorphic, then $E_\rho$ is Stein.
    \item If $s$ is holomorphic, then $s^*K_{E_\rho/X}$ is ample.
\end{enumerate}
\end{theorem}

The proof uses the positivity statement for exceptional conormal bundles in Proposition~\ref{prop:exceptional-conormal-big}. 

\subsection{Steinness of Affine Bundles over K\"ahler Manifolds}\label{subsec:intro-affine}

Let $E\to X$ be a holomorphic vector bundle. We write $H^1(X,E)$ for $H^1(X,\mathcal O(E))$. A class $t\in H^1(X,E)$ determines a holomorphic affine bundle $E_t$ modelled on $E$; the precise convention is recalled in \S\ref{subsec:holomorphic-fibre-bundles}.

Affine line bundles provide an important one-dimensional model. We shall use the following Stein criterion.

\begin{theorem}[Ohsawa {\normalfont\cite[p.~2814 and Theorem~2.1]{Ohsawa2018}}]\label{thm:riemann-affine-known}
Let $R$ be a compact Riemann surface, let $L\to R$ be a holomorphic line bundle with $\deg L\leq0$, and let $t\in H^1(R,L)$. Then the associated affine line bundle $L_t$ is weakly $1$-complete, and
\[
 L_t\text{ is Stein}\quad\Longleftrightarrow\quad t\ne0.
\]
\end{theorem}

See also Suzuki \cite[Proposition~(2.2)(2), pp.~187--188]{Suzuki1975}, Ueda \cite{Ueda1983}, and Koike--Ueda \cite[Theorems~1.1--1.2]{KoikeUeda2025} for related constructions of pseudoconvex neighbourhoods and plurisubharmonic exhaustion functions.

We seek a higher-dimensional analogue of Ohsawa's theorem. For vector bundles, the natural curvature counterpart of the condition $\deg L\le0$ is Griffiths seminegativity. In higher dimension, however, one must impose an additional condition on the affine cohomology class: we require it to admit a parallel Dolbeault representative.

Let $(E,h)$ be a Hermitian holomorphic vector bundle and let $D_h=D'_h+\overline\partial$ be its Chern connection. For an $E$-valued $(0,1)$-form $a$, we say that $a$ is parallel if $D_ha=0$. If $a$ is parallel, then
\[
\omega_{h,a}:=-\sqrti\,h(a,a)
\]
is a closed semipositive $(1,1)$-form. The Dolbeault class $[a]\in H^1(X,E)$ determines the affine bundle $E_{[a]}$.

\begin{theorem}\label{thm:main-affine}
Let $X$ be a weakly $1$-complete complex manifold, let $(E,h)$ be a Griffiths seminegative Hermitian holomorphic vector bundle over $X$, and let $a$ be a parallel $E$-valued $(0,1)$-form. Then $E_{[a]}$ is weakly $1$-complete.
If, in addition, $X$ is compact K\"ahler, then the following conditions are equivalent:
\begin{enumerate}
    \item $E_{[a]}$ is Stein;
    \item $[\omega_{h,a}]$ is a K\"ahler class.
\end{enumerate}
\end{theorem}

Local solutions of the Dolbeault equation define a smooth section $s$ of $E_{[a]}$; see Section~\ref{subsec:holomorphic-fibre-bundles}. The squared fibrewise $h$-distance from $s$ defines a function $\Phi$ on $E_{[a]}$. The main step in the proof is the computation of $\ddbar\,\Phi$.

When the base is a Riemann surface, the condition that $a$ be parallel is automatic for harmonic representatives. We therefore obtain the following higher-rank extension of Theorem~\ref{thm:riemann-affine-known}.

\begin{corollary}\label{cor:riemann-affine}
Let $X$ be a compact Riemann surface and let $(E,h)$ be a Griffiths seminegative Hermitian holomorphic vector bundle. For every $t\in H^1(X,E)$, the affine bundle $E_t$ is weakly $1$-complete. Moreover, the following conditions are equivalent:
\begin{enumerate}
    \item $t\ne0$;
    \item $E_t$ is Stein.
\end{enumerate}
\end{corollary}

The Bochner--Kodaira--Nakano identity gives another useful consequence. We check in Section~\ref{subsec:affine-proof} that, if $(E^*,h^*)$ is Nakano semipositive, then the harmonic representative of every class $t\in H^1(X,E)$ is parallel. Since Nakano semipositivity of $E^*$ implies Griffiths seminegativity of $E$, Theorem~\ref{thm:main-affine} yields the following.

\begin{corollary}\label{cor:nakano-affine-weak}
Let $X$ be a compact K\"ahler manifold and let $(E,h)$ be a Hermitian holomorphic vector bundle such that $(E^*,h^*)$ is Nakano semipositive. Then, for every $t\in H^1(X,E)$, the affine bundle $E_t$ is weakly $1$-complete.
\end{corollary}

Theorem~\ref{thm:main-affine} also recovers the theorem of Greb--Wong on canonical extensions: if $(X,\omega)$ is compact K\"ahler with nonnegative holomorphic bisectional curvature, then the canonical extension $\Omega^1_{[\omega]}$ is Stein \cite[Theorem~3.4]{GrebWong2020}. The term \emph{canonical extension} is sometimes used for extensions of general vector bundles. In this paper, following Greb--Wong, we use this term for the affine bundle $\Omega^1_{[\omega]}$ associated with a K\"ahler class $[\omega]$. Their proof uses the structure theorem for compact K\"ahler manifolds with nonnegative holomorphic bisectional curvature, and hence ultimately Mok's uniformisation theorem \cite[Main Theorem, p.~179]{Mok1988}. In contrast, Theorem~\ref{thm:main-affine} is based on an explicit construction of a plurisubharmonic exhaustion. One advantage of this proof is that, when combined with the heat-flow exhaustion estimate of Ni--Tam, it extends the result to complete K\"ahler manifolds.

\begin{theorem}\label{thm:main-complete}
Let $(X,\omega)$ be a complete K\"ahler manifold with nonnegative holomorphic bisectional curvature. Then the canonical extension $\Omega^1_{[\omega]}$ is Stein.
\end{theorem}

\paragraph*{Acknowledgements.}
The author is grateful to Hisamoto for his guidance and support during the preparation of this paper, and to Adachi, Lee, and Seo for helpful discussions and for reading the preprint. The author is supported by JST SPRING, Grant Number JPMJSP2156.

\paragraph*{Use of AI.}
During the preparation of this manuscript, the author used ChatGPT to improve the English exposition and to search the literature. The author takes full responsibility for the mathematical arguments and all other content of the paper.

\section{Preliminaries}\label{sec:preliminaries}

\subsection{Holomorphic fibre bundles}\label{subsec:holomorphic-fibre-bundles}

A holomorphic fibre bundle $p\colon Y\to X$ with fibre $F$ is a holomorphic map which is locally biholomorphic over $X$ to the projection $U\times F\to U$, with holomorphic transition maps in the fibre variable. In this subsection, we recall the definitions and basic properties of ball bundles and affine bundles.

\subsubsection{Ball bundles and the Toledo class}

Let $X$ be connected, let $\tX$ be its universal cover, and let
\[
 \rho\colon\pi_1(X)\longrightarrow\Aut(\B^m)=\PU(m,1)
\]
be a representation. The quotient
\begin{equation}\label{eq:ball-bundle}
 E_\rho:=\tX\times_\rho\B^m\longrightarrow X
\end{equation}
is a holomorphic fibre bundle with fibre $\B^m$. We write $[x,w]$ for the equivalence class of $(x,w)$ under
\[
 (x,w)\sim(\gamma x,\rho(\gamma)w),\qquad \gamma\in\pi_1(X).
\]

Let
\[
 H(Z,W):=\sum_{j=1}^{m}Z_j\overline{W_j}-Z_{m+1}\overline{W_{m+1}}
\]
be the standard Hermitian form of signature $(m,1)$ on $\C^{m+1}$. The projective action of $\PU(m,1)$ is induced by $H$-unitary linear transformations and therefore preserves the sign of $H(Z,Z)$. It defines
\begin{equation}\label{eq:projective-bundle}
 F_\rho:=\tX\times_\rho\Pp^m\longrightarrow X.
\end{equation}
In the projective model,
\[
 \B^m=\{[Z]\in\Pp^m;\ H(Z,Z)<0\},
\]
we have
\begin{equation}\label{eq:ball-as-domain}
 E_\rho=\{[x,Z]\in F_\rho;\ H(Z,Z)<0\}.
\end{equation}
Thus $E_\rho$ is a locally pseudoconvex domain in $F_\rho$.

For surface groups, Toledo associated to a representation the integral of the pullback of the invariant K\"ahler form by an equivariant map; see \cite[Introduction, p.~125]{Toledo1989}. The following class is its cohomological version, with the sign and normalization specified below.

\begin{definition}[Toledo class]\label{def:toledo-class}
Let $s\colon X\to E_\rho$ be a smooth section. Since $\B^m$ is contractible, any two sections are homotopic, and hence
\begin{equation}\label{eq:toledo-class}
 c_\rho:=c_1\bigl(s^*T_{E_\rho/X}\bigr)\in H^2(X,\R)
\end{equation}
is independent of $s$. Let $\omega_{\B}$ be the invariant K\"ahler form on $\B^m$, normalized so that the first Chern form of $T_{\B^m}$ is $-\omega_{\B}$:
\[
 \omega_{\B}=-\frac{m+1}{2\pi}\ddbar\log(1-|z|^2).
\]
In a local trivialization, the fibre component of $s$ defines $s^*\omega_{\B}$. These local forms agree on overlaps because $\omega_{\B}$ is invariant under holomorphic automorphisms of $\B^m$, so $s^*\omega_{\B}$ is well defined on $X$. $\omega_{\B}$ is a K\"ahler--Einstein metric with $\operatorname{Ric}(\omega_{\B})=-\omega_{\B}$ in the first-Chern-form normalization, and hence
\begin{equation}\label{eq:toledo-pullback}
 c_\rho=-[s^*\omega_{\B}].
\end{equation}
\end{definition}

\subsubsection{Harmonic and pluriharmonic sections}\label{subsec:harmonic-sections}

To relate the Toledo class to analytic properties of the total space, we use sections with suitable complex-analytic properties. For reductive monodromy, harmonic-map theory provides such distinguished sections.

Let $(M,g)$ be a K\"ahler manifold and $(N,h)$ a Riemannian manifold. For a smooth map $u\colon M\to N$, let $\nabla$ be the connection induced by the Levi--Civita connections of $g$ and $h$. Thus $\nabla du$ is a section of $T_M^*\otimes T_M^*\otimes u^*T_N$. We denote its trace with respect to $g$ by
\[
 \tau(u):=\operatorname{tr}_g\nabla du\in\Gamma(M,u^*T_N).
\]
The map $u$ is harmonic if $\tau(u)=0$, and pluriharmonic if, for every holomorphic map $f\colon\mathbb D\to M$, the composition $u\circ f$ is harmonic. Pluriharmonicity is independent of the choice of K\"ahler metric on $M$ and is a property of its complex structure.

A section $s\colon X\to E_\rho$ corresponds to a $\rho$-equivariant map
\[
 \widetilde s\colon\tX\longrightarrow\B^m,
 \qquad
 \widetilde s(\gamma x)=\rho(\gamma)\widetilde s(x).
\]
If $\omega_X$ is a K\"ahler metric on $X$ and $\pi\colon\tX\to X$ is the covering map, we call $s$ harmonic or pluriharmonic when $\widetilde s$ has the corresponding property with respect to $\pi^*\omega_X$ and the K\"ahler--Einstein metric $\omega_{\B}$ defined above.

A representation $\rho\colon\pi_1(X)\to\PU(m,1)$ is called reductive if the real Zariski closure of $\rho(\pi_1(X))$ is reductive. In the ball case, $\rho$ is non-reductive if and only if $\rho(\pi_1(X))$ has exactly one common fixed point on $\partial\B^m$.

\begin{proposition}[Corlette--Siu--Sampson {\normalfont\cite[Corollary~3.5]{Corlette1988}; \cite[Proposition~3, pp.~81--84]{Siu1980}; \cite[pp.~125--134]{Sampson1986}}]\label{prop:pluriharmonic-section}
Let $X$ be a compact K\"ahler manifold and let $\rho\colon\pi_1(X)\to\PU(m,1)$ be reductive. Then $E_\rho$ admits a pluriharmonic section.
\end{proposition}

For a pluriharmonic section $s$, set
\begin{equation}\label{eq:rank-section}
 \rank ds:=\max_{x\in\tX}\rank_{\R}d\widetilde s_x.
\end{equation}

When $\dim_{\C}X\ge2$, the following rigidity theorem for pluriharmonic sections plays an important role.

\begin{theorem}[Siu--Carlson--Toledo {\normalfont\cite[Proposition~3 and Theorem~5, pp.~81--89]{Siu1980}; \cite[Theorems~2.3 and~3.5, Corollary~3.7, pp.~180--183]{CarlsonToledo1989}; \cite[p.~519]{Toledo1999}}]\label{thm:rank-classification}
Let $X$ be a compact K\"ahler manifold and let $s\colon X\to E_\rho$ be a pluriharmonic section. Then:
\begin{enumerate}
 \item if $\rank ds=0$, then $s$ is locally constant;
 \item if $\rank ds=1$, the image of $\widetilde s$ is contained in a geodesic of $\B^m$;
 \item if $\rank ds=2$, then near every point where $d\widetilde s$ has real rank two there are a holomorphic submersion $f\colon U\to\mathbb D$ and a harmonic map $v\colon\mathbb D\to\B^m$ such that $\widetilde s|_U=v\circ f$;
 \item if $\rank ds>2$, $s$ is holomorphic or anti-holomorphic.
\end{enumerate}
\end{theorem}

Let $K_{\B}(z,w)$ denote the Bergman kernel of $\B^m$; thus
\[
 K_{\B}(z,w)=C_m(1-\langle w,z\rangle)^{-(m+1)}
\]
for a positive constant $C_m$. For a pluriharmonic section, consider
\begin{equation}\label{eq:ball-psi}
 \Psi(z,w)
 :=\frac{1}{2\pi}\log\frac{K_{\B}(z,z)K_{\B}(w,w)}{|K_{\B}(z,w)|^2}
 =\frac{m+1}{2\pi}\log\frac{|1-\langle w,z\rangle|^2}{(1-|w|^2)(1-|z|^2)}.
\end{equation}
It is invariant under the diagonal action of $\Aut(\B^m)$, and hence defines
\begin{equation}\label{eq:ball-Phi}
 \Phi_s([x,w]):=\Psi(\widetilde s(x),w)
\end{equation}
on $E_\rho$. For $\B^1$-bundles this construction is due to Diederich--Ohsawa \cite[Proposition~1 and the Key Lemma, pp.~823--829]{DiederichOhsawa1985}; Seo proved for bounded symmetric-domain bundles that $\Phi_s$ is a plurisubharmonic exhaustion \cite[Theorem~4.1]{Seo2022}.

When $s$ is anti-holomorphic, the Levi form of $\Phi_s$ takes a particularly simple form. To keep the base and fibre contributions visible, fix a flat local trivialization
\[
 p^{-1}(U)\simeq U\times\B^m
\]
and write a point as $(x,w)$, with $x\in U$ and $w\in\B^m$. If $\operatorname{pr}_{\B}\colon U\times\B^m\to\B^m$ is the second projection, we use the notation
\begin{equation}\label{eq:bergman-fibre-notation}
 \omega_{\B}(w):=\operatorname{pr}_{\B}^{*}\omega_{\B}.
\end{equation}
Since the transition maps of the flat bundle are locally constant elements of $\Aut(\B^m)$ and $\omega_{\B}$ is $\Aut(\B^m)$-invariant, these local forms agree on overlaps.

\begin{lemma}\label{lem:ball-hessian}
If $s\colon X\to E_\rho$ is anti-holomorphic, then, in every flat local trivialization as above,
\begin{equation}\label{eq:ball-hessian}
 \ddbar\,\Phi_s(x,w)
 =\omega_{\B}(w)-p^*\bigl(s^*\omega_{\B}\bigr)(x),
\end{equation}
where $\partial$ and $\overline\partial$ are taken on the total space $E_\rho$, in both the base and fibre variables.
\end{lemma}

\begin{proof}
In the local coordinates $(x,w)$, formula \eqref{eq:ball-psi} gives
\[
 \Phi_s(x,w)
 =\frac{m+1}{2\pi}\log|1-\langle w,s(x)\rangle|^2
 -\frac{m+1}{2\pi}\log(1-|w|^2)
 -\frac{m+1}{2\pi}\log(1-|s(x)|^2).
\]
Since $s$ is anti-holomorphic, $1-\langle w,s(x)\rangle$ is holomorphic in $(x,w)$, and hence
\[
 \ddbar\,\log|1-\langle w,s(x)\rangle|^2=0.
\]
The second term depends only on the fibre variable and satisfies
\[
 \ddbar\,\bigl(-\frac{m+1}{2\pi}\log(1-|w|^2)\bigr)=\omega_{\B}(w).
\]
The third term depends only on the base variable. Anti-holomorphicity of $s$ gives
\[
 \ddbar\,\bigl(-\frac{m+1}{2\pi}\log(1-|s(x)|^2)\bigr)
 =-p^*\bigl(s^*\omega_{\B}\bigr)(x).
\]
Adding the three terms proves \eqref{eq:ball-hessian}.
\end{proof}

\subsubsection{Affine bundles}

Let $E\to X$ be a holomorphic vector bundle of rank $r$. We write $H^1(X,E)$ for the sheaf cohomology $H^1(X,\cO(E))$. Let $\{U_i\}$ be a holomorphic trivializing cover with transition functions $g_{ij}$. A holomorphic affine bundle modelled on $E$ has local fibre coordinates
\begin{equation}\label{eq:affine-transition}
 w_j=g_{ij}w_i+t_{ij},
\end{equation}
where the holomorphic transition maps satisfy the cocycle condition and $t_{ij}\in\cO(U_i\cap U_j,E)$. Thus $\{t_{ij}\}$ defines a class $[t]\in H^1(X,E)$, and we denote the associated affine bundle by $E_t$.

Let $a\in A^{0,1}(X,E)$, $\overline\partial a=0$, represent $[t]$. Locally choose $f_i$ with $\overline\partial f_i=a|_{U_i}$. Then $t_{ij}:=f_j-g_{ij}f_i$ is holomorphic, and the $f_i$ define a smooth section $f\colon X\to E_{[a]}$ satisfying $\overline\partial f=a$.

\subsection{Stein and $1$-convex manifolds}\label{subsec:stein-preliminaries}

A complex manifold $Y$ is called Stein if it is holomorphically convex and holomorphically separable: the $\mathcal O(Y)$-convex hull of every compact subset is compact, and any two distinct points of $Y$ are separated by a function in $\mathcal O(Y)$. By Grauert's solution of the Levi problem, this is equivalent to $1$-completeness \cite[Theorem~2 and pp.~468--469]{Grauert1958}.

A complex manifold $Y$ is called $1$-convex if it admits a smooth exhaustion function which is strictly plurisubharmonic outside a compact set. Its Remmert reduction is a proper holomorphic map $r\colon Y\to S$ to a Stein space $S$ such that there is a finite set $\Sigma\subset S$ for which
\[
 r\colon Y\setminus A\longrightarrow S\setminus\Sigma,\qquad A:=r^{-1}(\Sigma),
\]
is biholomorphic \cite[\S2, especially Satz~5]{Grauert1962}. The compact analytic set $A$ is the exceptional set; it is empty exactly when $Y$ is Stein.

\begin{theorem}[Andreotti--Frankel--Hamm-type bound]\label{thm:oneconvex-topology}
Let $Y$ be a $1$-convex complex manifold of complex dimension $N$ and let $A$ be its exceptional set.
\begin{enumerate}
 \item If $A=\varnothing$, then $Y$ has the homotopy type of a CW complex of real dimension at most $N$.
 \item If $A\neq\varnothing$ and $d:=\dim_{\C}A$, then
 \[
 H_k(Y,\C)=0\qquad\text{for }k>\max\{N,2d\}.
 \]
\end{enumerate}
\end{theorem}

\begin{proof}
The first assertion is Andreotti--Frankel \cite[pp.~715--717]{AndreottiFrankel1959}. For the second, let $r\colon Y\to S$ be the Remmert reduction. Since $r$ contracts the connected components of $A$ to the finite set $\Sigma=r(A)$ and is a homeomorphism outside of $A$,
\[
 H_k(Y,A;\C)\simeq H_k(S,\Sigma;\C).
\]
Hamm's theorem gives $H_k(S,\C)=0$ for $k>N$ \cite[Satz~1]{Hamm1983}, while $H_k(A,\C)=0$ for $k>2d$. The long exact sequence of $(Y,A)$ gives the result.
\end{proof}

\begin{corollary}\label{cor:dimension-obstruction}
Let $X$ be a compact complex manifold of complex dimension $n$, and let $p\colon Y\to X$ be a holomorphic fibre bundle with contractible fibre of complex dimension $r$.
\begin{enumerate}
 \item If $Y$ is Stein, then $r\ge n$.
 \item If $Y$ is $1$-convex and $r<n$, then the exceptional set has an irreducible component $A$ of complex dimension at least $n$. If, in addition, every compact analytic subset of a fibre is zero-dimensional, then $\dim_{\C}A=n$ and $p|_A\colon A\to X$ is finite and surjective.
\end{enumerate}
\end{corollary}

\begin{proof}
Since the fibre is contractible, $p$ is a homotopy equivalence and 
\[
 H_{2n}(Y,\C)\simeq H_{2n}(X,\C)\ne0.
\] The theorem gives (1) and, when $r<n$, forces $2\dim_{\C}A\ge2n$. Since the fibres are contractible and, under the additional hypothesis in (2), contain no positive-dimensional compact analytic subsets, $p|_A$ has zero-dimensional fibres, hence is finite; surjectivity and $\dim_{\C}A=n$ follow.
\end{proof}

We next recall the holomorphic de Rham theorem. Define
\begin{equation}\label{eq:holomorphic-de-rham-cohomology}
 H^q_{\mathrm{hol}}(Y):=\frac{\ker\bigl(d:H^0(Y,\Omega_Y^q)\to H^0(Y,\Omega_Y^{q+1})\bigr)}{\operatorname{im}\bigl(d:H^0(Y,\Omega_Y^{q-1})\to H^0(Y,\Omega_Y^q)\bigr)}.
\end{equation}

\begin{theorem}[Serre {\normalfont\cite[\S2.1, pp.~213--218]{Serre1954}}]\label{thm:holomorphic-de-rham}
If $Y$ is Stein, then the natural map $H^q_{\mathrm{hol}}(Y)\to H^q_{\mathrm{dR}}(Y,\C)$ is an isomorphism for every $q\ge0$.
\end{theorem}

Finally, we record the positivity of conormal bundles of exceptional sets.  We use Ancona's notion of ampleness for coherent sheaves: if $Z$ is compact, a coherent sheaf $\mathcal S$ on $Z$ is called \emph{ample} if, for every coherent sheaf $\mathcal Q$ on $Z$, the sheaf $\mathcal Q\otimes\operatorname{Sym}^n\mathcal S$ is generated by global sections for all $n\gg0$; this is Definition~2.1 of \cite[p.~91]{Ancona1982}.  By \cite[Proposition~2.5, p.~92]{Ancona1982}, this is equivalent to cohomological ampleness.

\begin{theorem}[Ancona {\normalfont\cite[Theorem~3.7]{Ancona1982}}]\label{thm:ancona}
Let $Y$ be a $1$-convex complex space and let $A$ be its exceptional set. Then there exists a coherent ideal sheaf $\mathcal J\subset\mathcal O_Y$ defining $A'$ with $A'_{\mathrm{red}}=A$ such that $\mathcal J/\mathcal J^2$ is ample.
\end{theorem}

For a vector bundle $E$ we write $\Pp(E^*)$ for the bundle of one-dimensional subspaces of $E^*$ and $\mathcal O_{\Pp(E^*)}(1)$ for the dual tautological line bundle; with this convention, if $\pi:\Pp(E^*)\to Z$, then $\pi_*\mathcal O_{\Pp(E^*)}(m)=\operatorname{Sym}^mE$ for $m\ge0$.  If $Z$ is a compact irreducible complex space, we call $E$ \emph{big} when $\mathcal O_{\Pp(E^*)}(1)$ is big after pullback to a resolution.

\begin{proposition}[Exceptional conormal positivity]\label{prop:exceptional-conormal-big}
Let $Y$ be a $1$-convex complex space and let $A$ be a positive-dimensional irreducible component of its exceptional set. Assume that $A\subset Y$ is a local complete intersection, and write $N:=N_{A/Y}$. Then:
\begin{enumerate}
 \item $N^*$ is big, equivalently $\mathcal O_{\Pp(N)}(1)$ is big;
 \item if $\mathcal O_{\Pp(N)}(1)$ is semipositive, then $\det N^*$ is big.
\end{enumerate}
\end{proposition}

We believe that this result is well known to experts, but we have not found a suitable reference. The proof is therefore given in Appendix~\ref{app:exceptional-conormal}.

\subsection{The Demailly--P\u{a}un criterion}\label{subsec:demailly-paun}

The following K\"ahler criterion will be used in the proofs of Theorems~\ref{thm:main-ball} and~\ref{thm:main-affine}.

\begin{theorem}[Demailly--P\u{a}un {\normalfont\cite[Main Theorem~0.1]{DemaillyPaun2004}}]\label{thm:demailly-paun}
Let $X$ be a compact K\"ahler manifold. Set
\[
 \mathcal P:=\left\{\alpha\in H^{1,1}(X,\R);\ \int_Z\alpha^{\dim Z}>0\text{ for every irreducible analytic }Z\subset X,\ \dim Z>0\right\}.
\]
Then the K\"ahler cone of $X$ is one of the connected components of $\mathcal P$.
\end{theorem}

\begin{corollary}\label{cor:demailly-paun-semipositive}
Let $\alpha$ be a smooth closed semipositive $(1,1)$-form on a compact K\"ahler manifold $X$. Suppose that, for every positive-dimensional irreducible analytic subset $Z\subset X$, there exists $z\in Z_{\mathrm{reg}}$ such that $\alpha|_{T_zZ}$ is positive definite. Then $[\alpha]$ is a K\"ahler class.
\end{corollary}

\begin{proof}
If $q=\dim Z$, semipositivity and strict positivity at one smooth point imply $\int_Z\alpha^q>0$. Hence $[\alpha]\in\mathcal P$. Since $\alpha+t\omega_X$ is K\"ahler for every $t>0$, $[\alpha]$ lies in the K\"ahler component of $\mathcal P$.
\end{proof}

\section{Proofs of the main results}\label{sec:proof-main-results}

\subsection{Common section and rigidity lemmas}\label{subsec:section-lemma}

The following lemma is the common necessary condition for Steinness used for ball and affine bundles.

\begin{lemma}\label{lem:section-nondegeneracy}
Let $X$ be a compact K\"ahler manifold, let $p\colon Y\to X$ be a holomorphic fibre bundle, and let $s\colon X\to Y$ be a smooth section. Assume that $Y$ is Stein. Then, for every positive-dimensional compact irreducible analytic subset $Z\subset X$, there exists $z\in Z_{\mathrm{reg}}$ such that
\begin{equation}\label{eq:section-injectivity}
 (\overline\partial s)_z\colon T_{Z,z}^{0,1}\longrightarrow T_{Y/X,s(z)}
\end{equation}
is injective.
\end{lemma}

\begin{proof}
Let $q=\dim Z$, and let $\mu\colon\widehat Z\to Z$ be a resolution. Choose a K\"ahler form $\omega_X$ on $X$. Since $Y$ is Stein, Theorem~\ref{thm:holomorphic-de-rham} gives a closed holomorphic $2$-form $\xi\in H^0(Y,\Omega_Y^2)$ representing the de Rham class $[p^*\omega_X]$. Hence
\[
 [s^*\xi]=[\omega_X]\quad\text{in }H^2_{\mathrm{dR}}(X,\C),
\]
and therefore
\[
 \int_{\widehat Z}\mu^*(s^*\xi)^q
 =\int_{\widehat Z}\mu^*\omega_X^q>0.
\]
Thus the $(q,q)$-component of $(s^*\xi)^q$ is nonzero at some point $z_0\in Z_{\mathrm{reg}}$.

We now describe this component in local coordinates. Choose holomorphic coordinates $(z_1,\ldots,z_q)$ on $Z$ centred at $z_0$ and fibre coordinates $(w_1,\ldots,w_r)$ near $s(z_0)$, so that locally $s(z)=(z,s^1(z),\ldots,s^r(z))$. Along $Z$, write the holomorphic $2$-form $\xi$ in the form
\[
 \xi
 =\frac12\sum_{i,j}\xi_{ij}\,dz_i\wedge dz_j
 +\sum_{i,\alpha}\xi_{i\alpha}\,dz_i\wedge dw_\alpha
 +\frac12\sum_{\alpha,\beta}\xi_{\alpha\beta}\,dw_\alpha\wedge dw_\beta,
\]
where terms involving base coordinates transverse to $Z$ are omitted, since they vanish after restriction to $Z$. Since
\[
 s^*dw_\alpha
 =\sum_{j=1}^q\frac{\partial s^\alpha}{\partial z_j}\,dz_j
 +\sum_{j=1}^q\frac{\partial s^\alpha}{\partial\overline z_j}\,d\overline z_j,
\]
every $(0,1)$-factor occurring in $s^*\xi$ comes from the matrix
\[
 B=\left(\frac{\partial s^\alpha}{\partial\overline z_j}\right)_{\alpha,j},
\]
which represents $\overline\partial s|_{T_Z^{0,1}}$. Consequently the coefficient of
\[
 dz_1\wedge\cdots\wedge dz_q\wedge d\overline z_1\wedge\cdots\wedge d\overline z_q
\]
in the $(q,q)$-part of $(s^*\xi)^q$ is a linear combination of the $q\times q$ minors of $B$ (with coefficients depending on $\xi$ and on the $(1,0)$-part of $ds$). Equivalently, this $(q,q)$-component factors through $\bigwedge^q(\overline\partial s)_{z_0}$. Since it is nonzero at $z_0$, we must have
\[
 \bigwedge^q(\overline\partial s)_{z_0}\ne0.
\]
Hence $(\overline\partial s)_{z_0}\colon T_{Z,z_0}^{0,1}\to T_{Y/X,s(z_0)}$ is injective.
\end{proof}

\begin{lemma}\label{lem:harmonic-section-uniqueness}
Let $X$ be a compact K\"ahler manifold and let $E_\rho\to X$ be a locally trivial holomorphic ball bundle. If $E_\rho$ admits a pluriharmonic section which is not holomorphic, then $E_\rho$ admits no holomorphic section.
\end{lemma}

\begin{proof}
Suppose that $s$ is a non-holomorphic pluriharmonic section and that $t$ is a holomorphic section. Their lifts $\widetilde s,\widetilde t\colon\widetilde X\to\B^m$ are $\rho$-equivariant harmonic maps. Pointwise geodesic interpolation gives a $\rho$-equivariant homotopy between them. Hartman's uniqueness argument applies to equivariant maps because the distance and energy densities descend to the compact quotient $X$. By Assertions~(H) and~(I) of \cite[pp.~675--676]{Hartman1967}, either both maps are constant or, if they are distinct, their images lie in a real geodesic and their mutual distance is constant.

A holomorphic map with image in a real geodesic is locally constant, so $\widetilde t$ is constant. Then $\widetilde s$ is contained in the same real geodesic at fixed distance from $\widetilde t$, hence is constant as well. Thus $s$ is holomorphic, a contradiction.
\end{proof}

\begin{lemma}\label{lem:finite-multisection-obstruction}
Let $X$ be a positive-dimensional compact K\"ahler manifold and let $E_\rho\to X$ be a locally trivial holomorphic ball bundle admitting a non-holomorphic pluriharmonic section. Then there is no irreducible compact analytic subset $A\subset E_\rho$ such that $p|_A\colon A\to X$ is finite and surjective.
\end{lemma}

\begin{proof}
Assume that $A\subset E_\rho$ is an irreducible compact analytic subset such that $p|_A\colon A\to X$ is finite and surjective. Let $X'$ be the normalization of $A$. Then $q\colon X'\to X$ is a finite ramified covering, and the inclusion $A\subset E_\rho$ induces a holomorphic section of $q^*E_\rho$. The space $X'$ is a compact K\"ahler space by \cite[pp.~154--158]{AbeJinShima2004}; after resolving its singularities we obtain a compact K\"ahler manifold $\widehat X$ and a holomorphic section of the pulled-back ball bundle over $\widehat X$.

The pullback of the original pluriharmonic section remains non-holomorphic: the map $\widehat X\to X$ is generically finite, so holomorphicity of the pullback would imply holomorphicity of the original section on a dense open set and hence everywhere. This contradicts Lemma~\ref{lem:harmonic-section-uniqueness}.
\end{proof}

\subsection{Ball bundles}\label{subsec:ball-proof}

\subsubsection{Ball bundles over compact Riemann surfaces}\label{subsubsec:riemann-ball-proof}

We prove Theorem~\ref{thm:riemann-ball-criterion} by treating reductive and non-reductive monodromy separately. 
\paragraph{Reductive case over a Riemann surface.}
\mbox{}\par\nobreak\noindent
Assume that $\rho$ is reductive. Proposition~\ref{prop:pluriharmonic-section} provides a pluriharmonic section $s$. Since the base is a Riemann surface, pluriharmonicity is equivalent to harmonicity. We set
\[
 \operatorname{rk}ds:=\max_{x\in\widetilde R}\rank_{\R}d\widetilde s_x\in\{0,1,2\}.
\]
\begin{lemma}\label{lem:riemann-rankone}
If $E_\rho\to R$ admits a nonconstant harmonic section of real rank one, then $E_\rho$ is Stein.
\end{lemma}

\begin{proof}
By Theorem~\ref{thm:rank-classification}, the image of $\widetilde s$ is contained in a real geodesic $\ell\subset\B^m$. Since $s$ is nonconstant, $\ell$ and the unique complex geodesic $\Delta_\ell\simeq\mathbb D$ containing it are preserved by $\rho(\pi_1(R))$. Hence $D:=\widetilde R\times_\rho\Delta_\ell\subset E_\rho$ is a holomorphic $\B^1$-subbundle. The rank-one harmonic section is non-holomorphic, so Lemma~\ref{lem:harmonic-section-uniqueness} shows that $D$ admits no holomorphic section. Hence $D$ is Stein by Barrett's criterion \cite[Theorem~2]{Barrett1992}.

If $m=1$ we are done. For $m\ge2$, choose coordinates $(\zeta,\eta)\in\C\times\C^{m-1}$ with $\ell=(-1,1)\times\{0\}$ and $\Delta_\ell=\mathbb D\times\{0\}$. An automorphism preserving $\ell$ has the form
\begin{equation}\label{eq:riemann-geodesic-action}
 (\zeta,\eta)\longmapsto
 \left(f_\gamma(\zeta),\frac{U_\gamma\eta}{b_\gamma\zeta+a_\gamma}\right),
 \qquad
 f_\gamma(\zeta)=\varepsilon_\gamma\frac{a_\gamma\zeta+b_\gamma}{b_\gamma\zeta+a_\gamma},
\end{equation}
where $\varepsilon_\gamma\in\{\pm1\}$, $a_\gamma,b_\gamma\in\R$, $a_\gamma>0$, $a_\gamma^2-b_\gamma^2=1$, and $U_\gamma\in U(m-1)$. These transformations extend holomorphically from the ball to $\mathbb D\times\C^{m-1}$. Hence
\[
 F:=\bigl(\widetilde R\times\mathbb D\times\C^{m-1}\bigr)/\pi_1(R)\longrightarrow D
\]
is a holomorphic vector bundle. Since $D$ is Stein, so is $F$. The domain $E_\rho\subset F$ is locally pseudoconvex, hence Stein by Docquier--Grauert \cite[pp.~94--123]{DocquierGrauert1960}.
\end{proof}

For the real-rank-two case, the following result is crucial.

\begin{theorem}[Adachi--Lee--Seo {\normalfont\cite[Theorem~1.2]{AdachiLeeSeo2026}}]\label{thm:als-riemann}
Let $R$ be a compact Riemann surface of genus at least two. If a holomorphic ball bundle $E_\rho\to R$ admits a harmonic section having a point of real rank two, then $E_\rho$ is $1$-convex.
\end{theorem}

\begin{remark}\label{rem:riemann-ranktwo-genus}
The genus assumption in Theorem~\ref{thm:als-riemann} is automatic when the harmonic section has real rank two. For genus zero the monodromy is trivial, and the Bochner formula on the sphere forces every harmonic map to $\B^m$ to be constant.

Suppose now that the genus is one. Since harmonicity of maps from a Riemann surface is conformally invariant, we equip the base $R$ with a flat K\"ahler metric and denote its K\"ahler form by $\omega_R$. Let
\[
 \widetilde s:\C\longrightarrow\B^m
\]
be the equivariant harmonic map corresponding to the section. Because the monodromy acts by isometries of $\B^m$, the function
\[
 f:=|d\widetilde s|^2
\]
is invariant under the deck group and hence descends to a smooth function on $R$. Choose a local real orthonormal frame $e_1,e_2$ for the flat metric. Let $R^{\B}$ be the curvature tensor of the target, with the convention
\[
 \langle R^{\B}(u,v)v,u\rangle
 =K(u,v)\bigl(|u|^2|v|^2-\langle u,v\rangle^2\bigr).
\]
Here $K(u,v)<0$ is the sectional curvature when $u,v$ are linearly independent. The Bochner formula gives
\[
 \frac12\Delta f
 =|\nabla d\widetilde s|^2
 -\sum_{i,j=1}^{2}\langle R^{\B}(s_i,s_j)s_j,s_i\rangle,
\]
where $s_i(x):=d\widetilde s_x(e_i(x))\in T_{\widetilde s(x)}\B^m$ for $i=1,2$.
The right-hand side is nonnegative and is strictly positive wherever $d\widetilde s$ has real rank two. Thus $f$ is subharmonic on the compact Riemann surface $R$. By the maximum principle, $f$ is constant, so $\frac12\Delta f=0$. The Bochner formula therefore implies that $\operatorname{rk}ds\le1$. See \cite[\S3, formula~(16) and the corollary, pp.~123--124]{EellsSampson1964} for the Bochner argument.
\end{remark}

\begin{proposition}\label{prop:riemann-reductive}
Let $\rho\colon\pi_1(R)\to\PU(m,1)$ be reductive. Then $E_\rho$ is weakly $1$-complete. Moreover, if $E_\rho$ admits no holomorphic section, then it is Stein.
\end{proposition}

\begin{proof}
Weak $1$-completeness follows from Seo \cite[Theorem~4.1]{Seo2022}. Let $s$ be a harmonic section. Rank zero gives a holomorphic section and is excluded. The real-rank-one case follows from Lemma~\ref{lem:riemann-rankone}. If $\rank_{\R}ds=2$, Remark~\ref{rem:riemann-ranktwo-genus} and Theorem~\ref{thm:als-riemann} imply that $E_\rho$ is $1$-convex. Since $s$ is non-holomorphic, Lemma~\ref{lem:finite-multisection-obstruction} excludes every positive-dimensional compact analytic subset of $E_\rho$: any such subset maps finitely and surjectively to $R$. Hence the exceptional set is empty and $E_\rho$ is Stein.
\end{proof}

\paragraph{Non-reductive case over a Riemann surface.}
Assume now that $\rho$ is non-reductive. By \cite[Proposition~5.3]{Seo2022}, its image fixes a point of $\partial\B^m$. Sending this point to infinity by a Cayley transform, write the ball as the Siegel domain
\begin{equation}\label{eq:riemann-siegel}
 \mathcal H_m
 =\{(z,w)\in\C^{m-1}\times\C:\delta(z,w)>0\},
 \qquad
 \delta(z,w)=\operatorname{Im}w-|z|^2.
\end{equation}
The stabilizer of infinity acts by
\begin{equation}\label{eq:riemann-parabolic-action}
\begin{aligned}
 z&\longmapsto\lambda_\gamma U_\gamma z+a_\gamma,\\
 w&\longmapsto\lambda_\gamma^2w
 +2\sqrti\lambda_\gamma\langle U_\gamma z,a_\gamma\rangle
 +t_\gamma+\sqrti|a_\gamma|^2,
\end{aligned}
\end{equation}
where $\lambda_\gamma>0$, $U_\gamma\in U(m-1)$, $a_\gamma\in\C^{m-1}$, and $t_\gamma\in\R$; see \cite[Theorem~2.4]{Seo2022} for details. In particular,
\begin{equation}\label{eq:riemann-delta-transform}
 \delta(\rho(\gamma)(z,w))=\lambda_\gamma^2\delta(z,w).
\end{equation}
The transformations \eqref{eq:riemann-parabolic-action} extend affinely to $\C^m$. Let $W\to R$ be the affine bundle obtained from \eqref{eq:riemann-parabolic-action}, and let $A\to R$ be the affine bundle defined by its $z$-component. Then
\begin{equation}\label{eq:riemann-affine-tower}
 E_\rho\subset W\longrightarrow A\longrightarrow R.
\end{equation}
We separate the cases according to the dilation character $\lambda$.

\begin{proposition}\label{lem:riemann-lambda-one}
If $\rho$ is non-reductive and $\lambda\equiv1$, then $E_\rho$ is Stein.
\end{proposition}

\begin{proof}
If $A\to R$ has no holomorphic section, then, because $\lambda\equiv1$, its linear part is the unitary flat bundle associated with $U$. Corollary~\ref{cor:riemann-affine-proof} below gives that $A$ is Stein. The affine line bundle $W\to A$ is then holomorphically trivial as a torsor by Cartan's Theorem~B, so $W$ is Stein; since $E_\rho\subset W$ is locally pseudoconvex, Docquier--Grauert implies that $E_\rho$ is Stein. We may therefore assume that $A\to R$ has a holomorphic section $f$. In unitary flat coordinates write $f_i=U_{ij}f_j+a_{ij}$. Restrict $W\to A$ to $f(R)$. This gives an affine line bundle $T\to R$ with trivial linear part; write its transition functions as $w_i=w_j+c_{ij}(x)$. From \eqref{eq:riemann-parabolic-action} one obtains
\begin{equation}\label{eq:riemann-im-transition}
 \operatorname{Im}c_{ij}=|f_i|^2-|f_j|^2.
\end{equation}
Choose a smooth affine section of $T$ whose $\overline\partial$-differential is harmonic; in flat local coordinates its representatives $b_i$ are harmonic. Then $v:=\operatorname{Im}b_i-|f_i|^2$ is a globally defined function on $R$, and \eqref{eq:riemann-im-transition} gives
\begin{equation}\label{eq:riemann-flat-section}
 \ddbar\,v=-\sqrti\sum_\nu df_i^\nu\wedge d\overline{f_i^\nu}\le0.
\end{equation}
Integrating over $R$ gives $df_i=0$. Thus $f$ is flat. After a constant Heisenberg translation we may arrange $f=0$, and the monodromy has the form
\begin{equation}\label{eq:riemann-central-action}
 (z,w)\longmapsto(U_\gamma z,w+t_\gamma),
 \qquad t_\gamma\in\R.
\end{equation}
The character $t\colon\pi_1(R)\to\R$ is nonzero; otherwise the image would be conjugate into the compact subgroup $U(m-1)$ and hence reductive.

Let $T\to R$ be the affine line bundle with monodromy $w\mapsto w+t_\gamma$. It has no holomorphic section. Indeed, a holomorphic section would lift to a holomorphic function $g$ on $\widetilde R$ with $g(\gamma x)=g(x)+t_\gamma$. Then $dg$ would descend to a holomorphic one-form on $R$ all of whose periods are real. The imaginary part of $dg$ would be an exact harmonic one-form, hence zero; therefore $dg=0$ and $t=0$, a contradiction. By Ohsawa's affine-line criterion \cite[p.~2814]{Ohsawa2018}, $T$ is Stein. The projection $(z,w)\mapsto w$ identifies $W\to T$ with the holomorphic vector bundle associated with $U$, so $W$ is Stein. Since $E_\rho\subset W$ is locally pseudoconvex, Docquier--Grauert gives the conclusion.
\end{proof}

We next consider the case $\lambda\not\equiv1$. The remainder of this subsubsection
is devoted to the following proposition.

\begin{proposition}\label{lem:riemann-nontrivial-dilation}
If $\lambda\not\equiv1$, then $E_\rho$ is Stein.
\end{proposition}

Choose a harmonic real one-form on $R$ with periods $\log\lambda_\gamma$,
and let $u$ be its primitive on $\widetilde R$. Then
\begin{equation}\label{eq:riemann-u}
u_{x\bar x}=0,\qquad
u(\gamma x)=u(x)+\log\lambda_\gamma .
\end{equation}
Let $L$ be the flat line bundle associated with $\lambda$, and let $\mathcal U$
be the unitary flat bundle associated with $U$. The linear part $G=L\otimes\mathcal U$ of $A\to R$ is holomorphically unitary flat:
locally, writing $u=\operatorname{Re} h$ with $h$ holomorphic, the gauge $k=e^{-h}$
gives unitary flat frames.

\begin{lemma}\label{lem:riemann-positive-r}
There exists a positive smooth function $r$ on $\widetilde R$ such that
\begin{equation}\label{eq:riemann-r-properties}
r(\gamma x)=\lambda_\gamma r(x),\qquad r_{x\bar x}>0.
\end{equation}
\end{lemma}
\begin{proof}
Fix a Hermitian metric on $R$, put $q_0=|du|^2$, and let
$c_0=(\int_R q_0\,dV)/(\int_R dV)>0$.
Solve $\Delta v=c_0-q_0$ on $R$, and lift $v$ to $\widetilde R$.
For sufficiently small $\varepsilon>0$, set $r=e^u(1+\varepsilon v)$.
This function is positive and has the required automorphy. Since $\Delta u=0$,
\[
e^{-u}\Delta r
 =(1-\varepsilon+\varepsilon v)q_0+\varepsilon c_0
       +2\varepsilon\langle du,dv\rangle.
\]
Choose $\varepsilon$ so that $1-\varepsilon+\varepsilon v\ge3/4$.
Young's inequality gives
\[
e^{-u}\Delta r
 \ge\tfrac12 q_0+\varepsilon c_0
                 -4\varepsilon^2\|dv\|_\infty^2>0.
\]
On a Riemann surface, this is equivalent to $r_{x\bar x}>0$.
\end{proof}

The function $r$ will be used to construct a strictly plurisubharmonic exhaustion of $E_\rho$. First, set
\begin{equation}\label{eq:riemann-P-def}
P=\frac{r^2}{\delta}.
\end{equation}
Choose a smooth affine section $\sigma$ of $A$ whose $\bar\partial$-differential
is harmonic with respect to the unitary flat metric on $G$. In the local coordinates
$y=kz$, its representatives satisfy $\sigma_{x\bar x}=0$. Set
\begin{equation}\label{eq:riemann-Q-def}
Q=|kz-\sigma(x)|^2.
\end{equation}
We use the same notation for its pullback to $E_\rho$.

Next choose a smooth section $b$ of $W\to A$ which is affine in $z$:
\begin{equation}\label{eq:riemann-b-affine}
b_i(x,z)=\ell_i(x)z+c_i(x),\qquad
w_i-b_i=\lambda_{ij}^2(w_j-b_j).
\end{equation}
Such a section is obtained by a partition of unity depending only on the base variable.
Put
\begin{equation}\label{eq:riemann-PQJ-def}
F=w-b,\qquad a=e^{-4u},\qquad
T=1+a|F|^2,\qquad J=\log T.
\end{equation}
The functions $P,Q,J$ descend to smooth nonnegative functions on $E_\rho$, with
$P>0$, and $P+Q+J$ is an exhaustion. However, this sum need not be plurisubharmonic.
We shall therefore replace $P+Q$ by a suitable increasing convex function of $P+Q$.

Choose local coordinates $(U,x)$ on $R$ and a Siegel trivialization that extend
to a neighbourhood of $\overline U$. In these coordinates, put
\[
\theta_P=\frac{dw}{2i}-\sum_j\bar z_j\,dz_j-\delta\frac{r_x}{r}\,dx,
\qquad
\theta_{Q,j}=k\,dz_j+(k_xz_j-\sigma_{j,x})\,dx.
\]
Since
\[
\partial\delta=\frac{dw}{2i}-\sum_j\bar z_j\,dz_j,
\qquad \ddbar\delta=-\sum_j i\,dz_j\wedge d\bar z_j,
\]
direct differentiation gives
\[
\begin{aligned}
\ddbar P
&=\frac{2r}{\delta}\ddbar r
 +\frac2\delta i\,\partial r\wedge\bar\partial r
 -\frac{r^2}{\delta^2}\ddbar\delta\\
&\quad-\frac{2r}{\delta^2}i\bigl(
 \partial r\wedge\bar\partial\delta
 +\partial\delta\wedge\bar\partial r\bigr)
 +\frac{2r^2}{\delta^3}i\,\partial\delta\wedge\bar\partial\delta.
\end{aligned}
\]
Completing the square yields
\begin{equation}\label{eq:riemann-P-levi-formula}
\ddbar P
 =\frac{2rr_{x\bar x}}{\delta}\,i\,dx\wedge d\bar x
 +\frac{r^2}{\delta^2}\sum_j i\,dz_j\wedge d\bar z_j
 +\frac{2r^2}{\delta^3}\,i\theta_P\wedge\bar\theta_P
\end{equation}
For $Q$, direct differentiation gives the sum of squares
\begin{equation}\label{eq:riemann-Q-levi-formula}
\begin{aligned}
\ddbar Q
 &=\sum_j i\bigl(k\,dz_j+(k_xz_j-\sigma_{j,x})\,dx\bigr)\\
 &\qquad\wedge\overline{\bigl(k\,dz_j+(k_xz_j-\sigma_{j,x})\,dx\bigr)}
       +|\sigma_{\bar x}|^2 i\,dx\wedge d\bar x.
\end{aligned}
\end{equation}
Thus $P$ is strictly plurisubharmonic and $Q$ is plurisubharmonic.

\begin{proposition}\label{prop:riemann-local-estimates}
There is a constant $C=C_U\ge1$, independent of the point of $p^{-1}(U)$, such that
\begin{equation}\label{eq:riemann-PQ-weighted}
\ddbar(P+Q)\ge\frac1C\left(
 \frac{i\,dx\wedge d\bar x}{\sqrt T}
 +\sum_j i\theta_{Q,j}\wedge\bar\theta_{Q,j}
 +\frac{i\,dw\wedge d\bar w}{T^{3/2}}\right)
\end{equation}
and
\begin{equation}\label{eq:riemann-J-relative}
\ddbar J\ge-C\sqrt{1+Q}\,\ddbar(P+Q).
\end{equation}
\end{proposition}

\begin{proof}
All constants below depend only on the fixed local data and are absorbed into
a single $C$ at the end. Since $b$ is affine in $z$,
\[
\delta+|z|^2=\operatorname{Im}w\le |F|+C(1+|z|).
\]
Young's inequality and the definition of $Q$ therefore give
\begin{equation}\label{eq:riemann-growth}
\delta+|z|^2\le C\sqrt T,\qquad
|z|^2\le C(1+Q),\qquad (1+Q)^2\le CT.
\end{equation}
Completing squares in \eqref{eq:riemann-P-levi-formula}, using $r>0$ and
$r_{x\bar x}>0$, yields
\begin{equation}\label{eq:riemann-P-weighted}
\ddbar P\ge\frac1C\left(
 \frac{i\,dx\wedge d\bar x}{\delta}
 +\frac{\sum_j i\,dz_j\wedge d\bar z_j}{\delta^2}
 +\frac{i\,dw\wedge d\bar w}{\delta^2(\delta+|z|^2)}\right).
\end{equation}
Adding the exact sum of squares \eqref{eq:riemann-Q-levi-formula} to
\eqref{eq:riemann-P-weighted}, and using \eqref{eq:riemann-growth}, proves
\eqref{eq:riemann-PQ-weighted}.

For $J=\log T$, direct differentiation, using $u_{x\bar x}=0$, gives
\begin{align}
J_{w\bar x}
&=\frac a{T^2}\bigl(-\overline{b_x}-4u_{\bar x}\bar F
                    +a\bar F^{\,2}b_{\bar x}\bigr),
\label{eq:riemann-J-wx}\\
J_{z_j\bar x}
&=-b_{z_j}J_{w\bar x}-\frac{a\bar F}{T}b_{z_j\bar x},
\label{eq:riemann-J-zx}\\
J_{x\bar x}
&=\frac a{T^2}\bigl(|b_x+4u_xF|^2+|b_{\bar x}|^2\bigr)\notag\\
&\quad-\frac{2a^2}{T^2}\operatorname{Re}
 \bigl(\bar F^{\,2}(b_x+4u_xF)b_{\bar x}\bigr)\notag\\
&\quad-\frac{2a}{T}\operatorname{Re}
 \bigl(\bar F(b_{x\bar x}-4u_xb_{\bar x})\bigr).
\label{eq:riemann-J-xx}
\end{align}
Since $b$ is holomorphic in $z$, the vertical part is a square:
\begin{equation}\label{eq:riemann-J-levi-formula}
\begin{aligned}
\ddbar J
&=\frac a{T^2}i\left(dw-\sum_jb_{z_j}\,dz_j\right)
       \wedge\overline{\left(dw-\sum_jb_{z_j}\,dz_j\right)}\\
&\quad+J_{x\bar x}\,i\,dx\wedge d\bar x\\
&\quad+2\operatorname{Re}\left[
 i\left(J_{w\bar x}\,dw+\sum_jJ_{z_j\bar x}\,dz_j\right)\wedge d\bar x
 \right].
\end{aligned}
\end{equation}
Here $\operatorname{Re}\Omega=(\Omega+\bar\Omega)/2$ for a complex differential
form. Since
\[
b_x,b_{\bar x},b_{x\bar x}=O\bigl(\sqrt{1+Q}\bigr),\qquad
b_{z_j},b_{z_j\bar x}=O(1),\qquad F=O(\sqrt T),
\]
equations \eqref{eq:riemann-J-wx}--\eqref{eq:riemann-J-xx} imply
\begin{equation}\label{eq:riemann-J-orders}
\begin{aligned}
J_{w\bar x}&=O\left(\frac{\sqrt{1+Q}}T\right),\qquad
J_{z_j\bar x}=O\left(\frac1{\sqrt T}\right),\\
J_{x\bar x}&=O\left(\frac{\sqrt{1+Q}}{\sqrt T}\right).
\end{aligned}
\end{equation}
In the mixed terms of \eqref{eq:riemann-J-levi-formula}, substitute
\[
dz_j=k^{-1}\theta_{Q,j}-k^{-1}(k_xz_j-\sigma_{j,x})\,dx.
\]
Since $k_xz_j-\sigma_{j,x}=O(\sqrt{1+Q})$, the coefficient of
$i\,dx\wedge d\bar x$ outside the semipositive square becomes
\[
J_{x\bar x}-2\operatorname{Re}\left(
 k^{-1}\sum_j(k_xz_j-\sigma_{j,x})J_{z_j\bar x}\right)
 =O\left(\frac{\sqrt{1+Q}}{\sqrt T}\right),
\]
and the coefficient of $i\theta_{Q,j}\wedge d\bar x$ is
$k^{-1}J_{z_j\bar x}=O(T^{-1/2})$.
Cauchy--Schwarz applied to these mixed terms of $\ddbar J$, with the diagonal
weights in \eqref{eq:riemann-PQ-weighted}, now gives
\eqref{eq:riemann-J-relative}.
\end{proof}

\begin{proof}[Proof of Proposition~\ref{lem:riemann-nontrivial-dilation}]
Choose a finite cover $\{U_j\}_{j=1}^N$ of $R$ by coordinate neighbourhoods
as above. Let $C_j$ be the constant in \eqref{eq:riemann-J-relative} on
$p^{-1}(U_j)$, and put
\[
C=\max_{1\le j\le N}C_j.
\]
Set
\begin{equation}\label{eq:riemann-chi}
\chi(t)=\tau(1+t)^2,\qquad \tau>\max\{1,C/2\},
\end{equation}
and define
\begin{equation}\label{eq:riemann-final-exhaustion}
\Psi=\chi(P+Q)+J.
\end{equation}
Since $\chi$ is increasing and convex on $[0,\infty)$ and $Q\le P+Q$,
Proposition~\ref{prop:riemann-local-estimates} gives
\begin{equation}\label{eq:riemann-strict-positive}
\begin{aligned}
\ddbar\Psi
&\ge\bigl(2\tau(1+P+Q)-C\sqrt{1+Q}\bigr)\ddbar(P+Q)\\
&\ge(2\tau-C)(1+P+Q)\ddbar(P+Q)>0.
\end{aligned}
\end{equation}
Thus $\Psi$ is a strictly plurisubharmonic exhaustion of $E_\rho$,
and $E_\rho$ is Stein.
\end{proof}

Propositions~\ref{lem:riemann-lambda-one} and~\ref{lem:riemann-nontrivial-dilation} give the following.

\begin{proposition}\label{prop:riemann-nonreductive}
If $\rho\colon\pi_1(R)\to\PU(m,1)$ is non-reductive, then $E_\rho$ is Stein.
\end{proposition}

\Needspace{12\baselineskip}
\begin{theorem}[= Theorem~\ref{thm:riemann-ball-criterion}]
Let $R$ be a connected compact Riemann surface and let $p\colon E\to R$ be a holomorphic fibre bundle with fibre $\B^m$, $m\ge1$. Then $E$ is weakly $1$-complete. Moreover, the following conditions are equivalent:
\begin{enumerate}
    \item $E$ is Stein;
    \item $p\colon E\to R$ admits no holomorphic section;
    \item $E$ contains no positive-dimensional compact analytic subset.
\end{enumerate}
\end{theorem}

\begin{proof}
Combine Propositions~\ref{prop:riemann-reductive} and~\ref{prop:riemann-nonreductive}. Conversely, the image of a holomorphic section is a compact analytic curve, so a Stein total space admits no holomorphic section. A Stein manifold contains no positive-dimensional compact analytic subset, and the absence of such subsets in particular excludes holomorphic sections.
\end{proof}

\subsubsection{Ball bundles over higher-dimensional bases}

We now assume that the base has complex dimension at least two. The compact-analytic-subset criterion in Theorem~\ref{thm:riemann-ball-criterion} does not extend to this setting.

\begin{example}\label{ex:torus-B1-bundle}
There exists a holomorphic $\B^1$-bundle with reductive monodromy over a compact K\"ahler surface which contains no positive-dimensional compact analytic subset but is not Stein.
\end{example}

\begin{proof}
Let $X=\C^2/\Lambda$ be a two-dimensional complex torus with $\operatorname{NS}(X)=0$; such tori exist, for instance by Voisin's construction \cite[\S1, Lemma~1 and Remark~3, pp.~330--331]{Voisin2004}. Then $X$ has no positive-dimensional proper analytic subset: a compact curve would be an effective divisor with nonzero first Chern class.

Fix generators $\gamma_1,\gamma_2,\gamma_3,\gamma_4$ of $\pi_1(X)=\Lambda\simeq\mathbb Z^4$, and choose an infinite-order hyperbolic element $g\in\Aut(\B^1)$. Define
\[
 \rho(\gamma_1)=g,\qquad
 \rho(\gamma_j)=\operatorname{id}\quad(j=2,3,4),
 \qquad
 E_\rho:=\C^2\times_\rho\B^1\longrightarrow X.
\]
Since $g$ is semisimple, the Zariski closure of the cyclic group generated by $g$ is reductive; hence $\rho$ is reductive. Since the fibre has complex dimension one, Corollary~\ref{cor:dimension-obstruction} shows that $E_\rho$ is not Stein.

Suppose that $A\subset E_\rho$ were a positive-dimensional irreducible compact analytic subset. Since $\B^1$ contains no positive-dimensional compact analytic subset and $X$ has no positive-dimensional proper analytic subset, the projection $A\to X$ is finite and surjective. Let $\nu\colon A^\nu\to A$ be the normalization. Pulling back to the universal cover of $X$ gives a finite surjective map
\[
 \widetilde q\colon\widetilde A:=A^\nu\times_X\C^2\longrightarrow\C^2
\]
and a holomorphic map $F\colon\widetilde A\to\B^1$ induced by the fibre coordinate of $\C^2\times\B^1$; by equivariance, its restriction to a connected component of $\widetilde A$ is nonconstant. This is impossible, since that component is a finite branched covering of $\C^2$.
\end{proof}

Nevertheless, in the reductive higher-dimensional setting Steinness is governed by the Toledo class rather than by the absence of compact analytic subsets.

\Needspace{12\baselineskip}
\begin{theorem}[= Theorem~\ref{thm:main-ball}]
Let $X$ be a compact K\"ahler manifold with $\dim_{\C}X\geq 2$, and let $\rho\colon\pi_1(X)\to\Aut(\B^m)$ be reductive. Then the following are equivalent:
\begin{enumerate}
    \item $E_\rho$ is Stein;
    \item $c_\rho$ is a K\"ahler class.
\end{enumerate}
\end{theorem}

\begin{proof}
Assume first that $E_\rho$ is Stein, and let $s$ be the pluriharmonic section given by Proposition~\ref{prop:pluriharmonic-section}. By Lemma~\ref{lem:section-nondegeneracy}, $\overline\partial s_x$ is injective at some $x\in X$. Since $\dim_{\C}X\ge2$, $\rank_{\C}\overline\partial s_x\ge2$. The rank-zero and rank-one cases are impossible, while in the rank-two case Theorem~\ref{thm:rank-classification} locally factors $s$ through a Riemann surface and hence forces $\rank_{\C}\overline\partial s\le1$. Thus $\rank ds>2$. The same theorem shows that $s$ is holomorphic or anti-holomorphic, and the injectivity of $\overline\partial s_x$ excludes the holomorphic case. Hence $s$ is anti-holomorphic. In local fibre coordinates, $ds=\overline\partial s$.

Set
\begin{equation}\label{eq:ball-alpha}
 \alpha:=-s^*\omega_{\B}.
\end{equation}
Because $s$ is anti-holomorphic, $\alpha$ is semipositive. Let $Z\subset X$ be any positive-dimensional irreducible analytic subset. Lemma~\ref{lem:section-nondegeneracy} gives $z\in Z_{\mathrm{reg}}$ such that $\overline\partial s_z|_{T_z^{0,1}Z}$ is injective. At this point $\alpha|_{T_zZ}$ is positive definite. Corollary~\ref{cor:demailly-paun-semipositive}, together with \eqref{eq:toledo-pullback}, gives that
\[
 c_\rho=[\alpha]
\]
is a K\"ahler class.

Conversely, assume that $c_\rho$ is a K\"ahler class, and let $s$ be the pluriharmonic section given by Proposition~\ref{prop:pluriharmonic-section}. If $\rank ds\le2$, then
\[
 (s^*\omega_{\B})^n=0,
 \qquad n:=\dim_{\C}X\ge2.
\]
This contradicts
\[
 0<\int_Xc_\rho^n
 =\int_X(-s^*\omega_{\B})^n.
\]
Thus $\rank ds>2$. By Theorem~\ref{thm:rank-classification}, $s$ is holomorphic or anti-holomorphic. If $s$ were holomorphic, then $s^*\omega_{\B}$ would be semipositive. For any K\"ahler form $\omega_X$ on $X$,
\[
 \int_Xc_\rho\wedge\omega_X^{n-1}
 =-\int_Xs^*\omega_{\B}\wedge\omega_X^{n-1}\le0,
\]
which is impossible for the K\"ahler class $c_\rho$. Hence $s$ is anti-holomorphic.

By \eqref{eq:toledo-pullback} and the $\partial\overline\partial$-lemma, there exists a smooth real-valued function $\varphi$ on $X$ such that
\begin{equation}\label{eq:ball-base-correction}
 -s^*\omega_{\B}+\ddbar\,\varphi>0.
\end{equation}
By Lemma~\ref{lem:ball-hessian}, in a flat local trivialization with coordinates $(x,w)$,
\[
 \ddbar\,(\Phi_s+p^*\varphi)(x,w)
 =\omega_{\B}(w)
 +p^*\bigl(-s^*\omega_{\B}+\ddbar\,\varphi\bigr)(x)>0.
\]
Here the first term is positive definite in the fibre directions, while the second is positive definite in the base directions by \eqref{eq:ball-base-correction}.
The function $\Phi_s$ is exhaustive in the fibre direction, and $X$ is compact. Thus $\Phi_s+p^*\varphi$ is a strictly plurisubharmonic exhaustion function on $E_\rho$. Therefore $E_\rho$ is Stein.
\end{proof}

\Needspace{4\baselineskip}
\begin{example}\label{ex:stein-ball-bundles}
\leavevmode
\begin{enumerate}
    \item Let $X=\B^n/\Gamma$ be a compact ball quotient, $n\geq2$, where $\Gamma\subset\PU(n,1)$ is torsion-free, and let $c(z)=\overline z$. Define
    \[
     \rho\colon\Gamma\longrightarrow\PU(n,1),\qquad
     \rho(\gamma)=c\circ\gamma\circ c.
    \]
    The map $c\colon\B^n\to\B^n$ is $\rho$-equivariant, anti-holomorphic, and harmonic. Let $\omega_X$ be the standard K\"ahler--Einstein form on $X$, with the normalization used for $\omega_{\B}$. Since $c^*\omega_{\B}=-\omega_{\B}$, we have
    \[
     c_\rho=[\omega_X].
    \]
    Thus $E_\rho$ is Stein by Theorem~\ref{thm:main-ball}. More generally, if $p\colon Y\to X$ is a finite holomorphic branched covering or a holomorphic immersion, with $Y$ a compact complex manifold, then $p^*E_\rho\to Y$ is also a Stein ball bundle.

    \item A less elementary example is provided by the Mostow--Siu surfaces. Mostow and Siu constructed compact K\"ahler surfaces $X$ that are not ball quotients and admit K\"ahler metrics of negative sectional curvature \cite{MostowSiu1980}. Their construction gives a representation
    \[
     \rho\colon\pi_1(X)\longrightarrow\PU(2,1)
    \]
    with nondiscrete image and a $\rho$-equivariant holomorphic map $f\colon\tX\to\B^2$; see \cite[Theorem~2.4 and Proposition~2.6]{Deraux2005}. The monodromy image is Zariski dense, and hence $\rho$ is reductive; see \cite[Lemma~(11.5), p.~66]{DeligneMostow1986}.

    Let $\alpha$ denote the smooth closed $(1,1)$-form on $X$ descended from $f^*\omega_{\B}$. Its class is K\"ahler by Corollary~\ref{cor:demailly-paun-semipositive} and the local description in \cite[Proposition~2.6, pp.~8--9]{Deraux2005}. Define $\rho^*(\gamma)=c\circ\rho(\gamma)\circ c$. Then $c\circ f$ is $\rho^*$-equivariant and anti-holomorphic, and
    \[
     c_{\rho^*}=[\alpha].
    \]
    Theorem~\ref{thm:main-ball} therefore shows that $E_{\rho^*}$ is Stein. Since its monodromy image is nondiscrete, this bundle cannot be obtained by pulling back any of the standard ball bundles in (1). Indeed, those bundles and their pullbacks have discrete monodromy, and a holomorphic isomorphism of ball bundles conjugates their monodromy representations.

    For comparison, consider the classical picture for negative line bundles. If the total space of a line bundle $A\to X$ is $1$-convex with exceptional zero section, then $A^*$ is ample. By Kodaira's theorem, for sufficiently large $k$ there is a holomorphic immersion $j\colon X\longrightarrow\Pp^N$ such that
    \[
     A^{\otimes k}\simeq j^*\cO_{\Pp^N}(-1).
    \]
    In this sense, projective space, equipped with its tautological line bundle, is a classifying space for negative line bundles after taking tensor powers. It is natural to seek an analogous classifying space carrying a universal ball bundle for Stein ball bundles. The preceding example shows that ball quotients equipped with their standard ball bundles do not suffice, and finding such a classifying space appears difficult.
\end{enumerate}
\end{example}

The dimension assumption in Theorem~\ref{thm:main-ball} is essential, as the following example shows.

\begin{example}\label{ex:curve-rankone-zero-toledo}
There exists a locally trivial holomorphic $\B^1$-bundle over a compact Riemann surface with reductive monodromy and $c_\rho=0$ whose total space is Stein.
\end{example}

\begin{proof}
Let $R=\C/(\mathbb Z+\tau\mathbb Z)$, where $\operatorname{Im}\tau>0$, and fix the generators $\gamma_1(z)=z+1$ and $\gamma_2(z)=z+\tau$ of $\pi_1(R)\simeq\mathbb Z^2$. Choose a hyperbolic element $g\in\Aut(\B^1)$ and a unit-speed parametrization $c\colon\R\to\B^1$ of its invariant real geodesic such that $g(c(t))=c(t+\ell)$ for some $\ell>0$. Define
\[
 \rho(\gamma_1)=g,\qquad \rho(\gamma_2)=\operatorname{id}.
\]
The representation $\rho$ is reductive because $g$ is semisimple. The real-valued harmonic function
\[
 u(z)=\ell\left(\operatorname{Re}z-
 \frac{\operatorname{Re}\tau}{\operatorname{Im}\tau}\operatorname{Im}z\right)
\]
satisfies $u(z+1)=u(z)+\ell$ and $u(z+\tau)=u(z)$. Thus $\widetilde s:=c\circ u$ is a $\rho$-equivariant harmonic map of real rank one and defines a harmonic section $s$ of $E_\rho\to R$. Its image lies in a real geodesic, so $s^*\omega_{\B}=0$ and $c_\rho=0$. Nevertheless, $E_\rho$ is Stein by Lemma~\ref{lem:riemann-rankone}.
\end{proof}

We next prove the dimension obstruction for bounded symmetric-domain bundles.

\Needspace{12\baselineskip}
\begin{theorem}[= Theorem~\ref{thm:bsd-oneconvex-dimension}]
Let $X$ be a (not necessarily K\"ahler) compact complex manifold of dimension $n$, let $\Omega$ be a bounded symmetric domain of dimension $m$, and let $p\colon E\to X$ be a holomorphic fibre bundle with fibre $\Omega$. If $E$ is $1$-convex, then $m\ge n$.
\end{theorem}

\begin{proof}
Assume, to the contrary, that $m<n$. Since a bounded symmetric domain is contractible, $p\colon E\to X$ is a homotopy equivalence. If $E$ were Stein, Corollary~\ref{cor:dimension-obstruction} would already give $m\ge n$. Hence $E$ is not Stein. The same corollary gives an irreducible component $A$ of the exceptional set with
\[
 \dim_{\C}A=n,
\]
and $p|_A\colon A\to X$ finite and surjective, because a bounded symmetric domain contains no positive-dimensional compact analytic subset.

Let $\nu\colon A^\nu\to A$ be the normalization and put $q=p\circ\nu\colon A^\nu\to X$. Consider the pull-back bundle
\[
 E^\nu:=q^*E\longrightarrow A^\nu.
\]
The natural map $E^\nu\to E$ is finite. Hence $E^\nu$ is $1$-convex: indeed, after composing with the Remmert reduction of $E$, its Stein factorization has Stein target. The inclusion $A\subset E$ induces a tautological holomorphic section
\[
 s\colon A^\nu\longrightarrow E^\nu.
\]
Its image $S:=s(A^\nu)$ is a positive-dimensional irreducible component of the exceptional set of $E^\nu$. To see this, note that $S$ is compact analytic and hence contained in the exceptional set. Any irreducible exceptional component containing $S$ has finite fibres over $A^\nu$, since its fibres are compact analytic subsets of $\Omega$, and therefore has dimension at most $n=\dim S$. Thus it coincides with $S$. Moreover, as the graph of a section, $S\subset E^\nu$ is a local complete intersection.

Let $N:=N_{S/E^\nu}\simeq s^*T_{E^\nu/A^\nu}$. Equip $\Omega$ with its Bergman metric and normalize its invariant K\"ahler form $\omega_\Omega$ so that the induced Chern form of the canonical bundle is $\omega_\Omega$. The Bergman metric is $\Aut(\Omega)$-invariant, K\"ahler--Einstein, and has nonpositive holomorphic bisectional curvature; equivalently, $T_\Omega$ is Griffiths seminegative \cite{Mok1986}. It therefore induces a Griffiths seminegative metric on $N$, so the tautological line bundle $\mathcal O_{\Pp(N)}(1)$ is semipositive. Proposition~\ref{prop:exceptional-conormal-big} gives $L:=\det N^*=s^*K_{E^\nu/A^\nu}$ big on $A^\nu$. The induced semipositive Chern form of $L$ is locally $f^*\omega_\Omega$, where $s$ is represented by a holomorphic map $f$ to $\Omega$; in particular, it has pointwise rank at most $m$.

Since $L$ is big, $A^\nu$ is Moishezon. Choose a projective modification $\mu\colon M\to A^\nu$ with $M$ smooth projective. Then $\mu^*L$ is big and nef. Therefore, writing $f^*\omega_\Omega$ for the globally defined form obtained from the local representatives of the section,
\[
 \int_M\bigl(\mu^*f^*\omega_\Omega\bigr)^n
 =\int_M c_1(\mu^*L)^n>0.
\]
But $f^*\omega_\Omega$ has rank at most $m<n$, so the left-hand side is zero, a contradiction. Hence $m\ge n$.
\end{proof}

\Needspace{12\baselineskip}
\begin{theorem}[= Theorem~\ref{thm:main-ball-oneconvex}]
Let $X$ be a compact K\"ahler manifold with $\dim_{\C}X\geq 2$, and let $\rho\colon\pi_1(X)\to\Aut(\B^m)$ be reductive. Suppose that $E_\rho$ is $1$-convex. Then $E_\rho$ admits a section $s\colon X\to E_\rho$ which is either holomorphic or anti-holomorphic. Moreover, one of the following holds:
\begin{enumerate}
    \item If $s$ is anti-holomorphic, then $E_\rho$ is Stein.
    \item If $s$ is holomorphic, then $s^*K_{E_\rho/X}$ is ample.
\end{enumerate}
\end{theorem}

\begin{proof}
Let $s$ be the pluriharmonic section supplied by Proposition~\ref{prop:pluriharmonic-section}. If $E_\rho$ is Stein, the proof of Theorem~\ref{thm:main-ball} shows that $s$ is anti-holomorphic. Suppose now that $E_\rho$ is not Stein, and let $A$ be a positive-dimensional irreducible component of its exceptional set. Theorem~\ref{thm:ancona} and \cite[Propositions~2.5 and~2.8, pp.~92--93]{Ancona1982} imply that $A$ is Moishezon. Choose an irreducible compact curve in $A$, and let $j\colon C\to A\subset E_\rho$ be its normalization. Put $q=p\circ j$ and $Z=q(C)$. Since a ball contains no compact analytic curve, $q\colon C\to Z$ is finite. Thus the pullback $F:=q^*E_\rho\to p^{-1}(Z)$ is finite, and $F$ is $1$-convex by the same Stein-factorization argument used above.

The map $j$ induces a holomorphic section $t\colon C\to F$. Its image is a smooth local complete intersection and an irreducible component of the exceptional set of $F$: every compact analytic subset of $F$ has finite fibres over $C$ and hence dimension at most one. The normal bundle $N_{t(C)/F}\simeq t^*T_{F/C}$ has the Griffiths seminegative metric induced by the Bergman metric. Proposition~\ref{prop:exceptional-conormal-big} therefore gives that $L:=t^*K_{F/C}$ is big, so
\[
 0<\deg L=\int_C t^*\omega_{\B}.
\]
In particular, $t$ is not locally constant in flat trivializations. The harmonic sections $q^*s$ and $t$ agree by Hartman's uniqueness argument \cite[Assertion~(H) and its corollary, p.~675]{Hartman1967}, since a nonconstant holomorphic map has real rank two somewhere and cannot have image in a real geodesic.

If $\rank ds>2$, Theorem~\ref{thm:rank-classification} shows that $s$ is holomorphic or anti-holomorphic; the nonconstant holomorphic pullback $q^*s=t$ excludes the latter. Otherwise, choose $c\in C$ where the fibre component of $t$ has nonzero differential. Near $q(c)$, Theorem~\ref{thm:rank-classification}(3) gives a local factorization $s=v\circ f$, in flat coordinates, with $f$ a holomorphic submersion to a disc. Since $f\circ q$ has nonzero differential at $c$, it is locally biholomorphic, and $v=t\circ(f\circ q)^{-1}$ is holomorphic. Thus $s$ is holomorphic on a nonempty open set, and hence on all of $X$. Since $q^*s=t$ is nonconstant, $s$ cannot also be anti-holomorphic. This proves the first assertion and the anti-holomorphic case of the theorem.

Assume now that $s$ is holomorphic. Put $N:=s^*T_{E_\rho/X}$ and $L:=\det N^*=s^*K_{E_\rho/X}$. The complex hyperbolic metric has Griffiths negative tangent curvature. Therefore $N$ is Griffiths seminegative and $N^*$ is Griffiths semipositive. Equivalently, the induced metric on $\mathcal O_{\Pp(N)}(1)$ has semipositive curvature; hence $N^*$ is nef, and so is $L=\det N^*$.

Let $Z\subset X$ be a positive-dimensional irreducible analytic subset, and put $E_Z:=p^{-1}(Z)$. The restriction of a $1$-convex exhaustion of $E_\rho$ shows that $E_Z$ is $1$-convex. Since $s(Z)$ is compact and analytic, it is contained in the exceptional set of $E_Z$. Moreover, every compact analytic subset of a fibre $\B^m$ is finite, so no irreducible component of the exceptional set of $E_Z$ containing $s(Z)$ can have dimension larger than $\dim Z$. Thus $s(Z)$ is an irreducible component of that exceptional set.

The section $s(Z)\subset E_Z$ is a local complete intersection, and its conormal bundle is canonically $N^*_{s(Z)/E_Z}\simeq N^*|_Z$. The induced metric on $\mathcal O_{\Pp(N|_Z)}(1)$ is semipositive. Proposition~\ref{prop:exceptional-conormal-big} therefore gives directly that $L|_Z=\det(N^*|_Z)$ is big. Let $\mu\colon\widehat Z\to Z$ be a resolution. The pullback $\mu^*(L|_Z)$ is big and nef. Since $\widehat Z$ is compact K\"ahler,
\[
 \int_Z c_1(L)^{\dim Z}
 =\int_{\widehat Z}c_1(\mu^*L)^{\dim Z}>0.
\]
This holds for every positive-dimensional irreducible analytic subset $Z\subset X$. Since $L$ is nef, the classes $c_1(L)+\varepsilon[\omega_X]$ are K\"ahler for every $\varepsilon>0$. Theorem~\ref{thm:demailly-paun} therefore implies that $c_1(L)$ itself is a K\"ahler class. Choose any Hermitian metric $h_0$ on $L$. By the $\partial\overline\partial$-lemma, the K\"ahler representative of $c_1(L)$ is the curvature form of a metric $h_0e^{-\varphi}$ for a suitable smooth real function $\varphi$. Thus $L$ is positive, and the Kodaira embedding theorem \cite[Main Theorem, pp.~28--48]{Kodaira1954} implies that $L=s^*K_{E_\rho/X}$ is ample.

\end{proof}

\subsection{Affine bundles}\label{subsec:affine-proof}

Let $(E,h)$ be a Hermitian holomorphic vector bundle over a complex manifold $X$, let $a\in A^{0,1}(X,E)$ satisfy $\overline\partial a=0$, and let $f\colon X\to E_{[a]}$ be the smooth section constructed in Section~\ref{subsec:holomorphic-fibre-bundles}. The squared fibrewise distance from $f$ is the function
\begin{equation}\label{eq:affine-distance}
 \Phi\colon E_{[a]}\longrightarrow\R_{\ge0},
 \qquad \Phi(x,w):=h_x\bigl(w-f(x),w-f(x)\bigr).
\end{equation}
Write $u:=w-f(x)$, a smooth section of $p^*E$, and let $\widetilde D_h$ be the Chern connection on $p^*E$ induced by $h$. In a local holomorphic frame,
\[
 \widetilde D'_h=\partial+p^*(h^{-1}\partial_Xh),
\]
where $\partial$ differentiates on the total space $E_{[a]}$ and $\partial_X$ differentiates on the base. Below, base forms and curvature are pulled back by $p$ without further notation.

\begin{lemma}\label{lem:affine-hessian}
The Levi form of \eqref{eq:affine-distance} is
\begin{equation}\label{eq:affine-hessian}
\begin{aligned}
 \ddbar\,\Phi
 ={}&-h\bigl(\sqrti\,\Theta_hu,u\bigr)
 +\sqrti\,h\bigl(\widetilde D'_hu,\widetilde D'_hu\bigr)\\
 &-\sqrti\,h(a,a)
 +2\operatorname{Im}h(\widetilde D'_ha,u).
\end{aligned}
\end{equation}
Here $\Theta_h=D_h^2\in A^{1,1}(X,\operatorname{End}E)$ is the Chern curvature. The Hermitian pairing of bundle-valued forms includes the wedge product of their form components.
\end{lemma}

\begin{proof}
The computation is local. Fix $x_0\in X$ and choose a holomorphic frame of $E$ which is normal for $h$ at $x_0$. Then, at $x_0$,
\[
 \overline\partial u=-a,
 \qquad
 \overline\partial(\widetilde D'_hu)=\Theta_hu+\widetilde D'_ha,
\]
and metric compatibility gives
\[
 \overline\partial\Phi
 =-h(a,u)+h(u,\widetilde D'_hu).
\]
Differentiating the two terms separately gives
\begin{align*}
 -\sqrti\,\partial h(a,u)
 &=-\sqrti\,h(\widetilde D'_ha,u)-\sqrti\,h(a,a),\\
 \sqrti\,\partial h(u,\widetilde D'_hu)
 &=-h\bigl(\sqrti\,\Theta_hu,u\bigr)
   +\sqrti\,h\bigl(\widetilde D'_hu,\widetilde D'_hu\bigr)
   +\sqrti\,\overline{h(\widetilde D'_ha,u)}.
\end{align*}
Adding these identities gives \eqref{eq:affine-hessian}.
\end{proof}

\begin{lemma}\label{lem:omega-ha}
Assume
\begin{equation}\label{eq:affine-parallel-assumption}
 \overline\partial a=0,
 \qquad
 D'_ha=0.
\end{equation}
Then the form
\begin{equation}\label{eq:omega-ha}
 \omega_{h,a}:=-\sqrti\,h(a,a)
\end{equation}
is closed and semipositive. If, in addition, $(E,h)$ is Griffiths seminegative, then
\begin{equation}\label{eq:affine-hessian-positive}
 \ddbar\,\Phi
 =-h\bigl(\sqrti\,\Theta_hu,u\bigr)
 +\sqrti\,h\bigl(\widetilde D'_hu,\widetilde D'_hu\bigr)
 +p^*\omega_{h,a}
 \ge p^*\omega_{h,a}\ge0.
\end{equation}
\end{lemma}

\begin{proof}
The two equations in \eqref{eq:affine-parallel-assumption} imply $D_ha=0$. Compatibility of the Chern connection with $h$ therefore gives $d\,h(a,a)=0$, so $\omega_{h,a}$ is closed. Locally write $a=\sum_j a_j\,d\overline z_j$. Then
\[
 h(a,a)=\sum_{j,k}h(a_j,a_k)\,d\overline z_j\wedge dz_k
       =-\sum_{j,k}h(a_j,a_k)\,dz_k\wedge d\overline z_j.
\]
Thus
\[
 \omega_{h,a}=\sqrti\sum_{j,k}h(a_j,a_k)\,dz_k\wedge d\overline z_j\ge0;
\]
indeed, $-\sqrti\,\omega_{h,a}(v,\overline v)=|a(\overline v)|_h^2$ for $v\in T_x^{1,0}X$. Under Griffiths seminegativity, the first term on the right-hand side of \eqref{eq:affine-hessian} is semipositive. The second term is also semipositive, and the mixed term vanishes because $D'_ha=0$. This gives \eqref{eq:affine-hessian-positive}.
\end{proof}

\Needspace{12\baselineskip}
\begin{theorem}[= Theorem~\ref{thm:main-affine}]
Let $X$ be a weakly $1$-complete complex manifold, let $(E,h)$ be a Griffiths seminegative Hermitian holomorphic vector bundle over $X$, and let $a$ be a parallel $E$-valued $(0,1)$-form. Then $E_{[a]}$ is weakly $1$-complete.
If, in addition, $X$ is compact K\"ahler, then the following conditions are equivalent:
\begin{enumerate}
    \item $E_{[a]}$ is Stein;
    \item $[\omega_{h,a}]$ is a K\"ahler class.
\end{enumerate}
\end{theorem}

\begin{proof}
By Lemma~\ref{lem:omega-ha}, $\Phi$ is plurisubharmonic. Let $\psi$ be a smooth plurisubharmonic exhaustion function on $X$. Then $\Phi+p^*\psi$ is a plurisubharmonic exhaustion function on $E_{[a]}$. Thus $E_{[a]}$ is weakly $1$-complete.

Assume now that $X$ is compact K\"ahler. Suppose first that $E_{[a]}$ is Stein. Let $Z\subset X$ be a positive-dimensional irreducible analytic subset. Lemma~\ref{lem:section-nondegeneracy}, applied to $f$, gives a point $z\in Z_{\mathrm{reg}}$ such that
\[
 \overline\partial f_z=a_z\colon T_z^{0,1}Z\longrightarrow E_z
\]
is injective. Hence $\omega_{h,a}|_{T_zZ}$ is positive definite. Corollary~\ref{cor:demailly-paun-semipositive} implies that $[\omega_{h,a}]$ is a K\"ahler class.

Conversely, assume that $[\omega_{h,a}]$ is a K\"ahler class. By the $\partial\overline\partial$-lemma, there exists a smooth real-valued function $\varphi$ on $X$ such that
\begin{equation}\label{eq:affine-base-correction}
 \omega_{h,a}+\ddbar\,\varphi>0.
\end{equation}
Lemma~\ref{lem:omega-ha} then gives
\[
\begin{aligned}
 \ddbar\,(\Phi+p^*\varphi)
 &\ge \sqrti\,h\bigl(\widetilde D'_hu,\widetilde D'_hu\bigr)\\
 &\quad+p^*\bigl(\omega_{h,a}+\ddbar\,\varphi\bigr)>0.
\end{aligned}
\]
The right-hand side is positive definite: the second term is positive on the base, and the first is positive in the vertical directions.
Since $X$ is compact, $\Phi+p^*\varphi$ is exhaustive. Therefore $E_{[a]}$ is Stein.
\end{proof}

\Needspace{10\baselineskip}
\begin{corollary}[= Corollary~\ref{cor:riemann-affine}]\label{cor:riemann-affine-proof}
Let $X$ be a compact Riemann surface and let $(E,h)$ be a Griffiths seminegative Hermitian holomorphic vector bundle. For every $t\in H^1(X,E)$, the affine bundle $E_t$ is weakly $1$-complete. Moreover, the following conditions are equivalent:
\begin{enumerate}
    \item $t\ne0$;
    \item $E_t$ is Stein.
\end{enumerate}
\end{corollary}

\begin{proof}
Fix a K\"ahler metric on $X$. Let $\delta'_h=(D'_h)^*$ and $\overline\partial^{\,*}$ be the formal adjoints of $D'_h$ and $\overline\partial$ with respect to this metric and $h$, and set
\[
 \Delta'=D'_h\delta'_h+\delta'_hD'_h,
 \qquad
 \Delta''=\overline\partial\,\overline\partial^{\,*}
             +\overline\partial^{\,*}\overline\partial.
\]
Let $a$ be the $\overline\partial$-harmonic representative of $t$. On the Riemann surface $X$, the Bochner--Kodaira--Nakano identity gives
\[
 0=\langle\Delta''a,a\rangle
 =\langle\Delta'a,a\rangle
  +\bigl\langle[\sqrti\Theta_h(E),\Lambda]a,a\bigr\rangle.
\]
For an $E$-valued $(0,1)$-form on a curve, $\Lambda a=0$ and $\Theta_h(E)\wedge a=0$, so the curvature commutator vanishes for bidegree reasons. Moreover $\delta'_ha=0$, again by bidegree. Hence
\[
 0=\|D'_ha\|^2+\|\delta'_ha\|^2=\|D'_ha\|^2.
\]
Thus $D'_ha=0$; together with $\overline\partial a=0$, this gives $D_ha=0$. Theorem~\ref{thm:main-affine} therefore applies, so $E_t$ is weakly $1$-complete. If $t\ne0$, then $a\not\equiv0$, so the semipositive form $\omega_{h,a}$ has positive integral and hence represents a K\"ahler class on the curve. Conversely, if $t=0$, then $E_t$ is isomorphic to its model vector bundle $E$. Its zero section is a compact complex curve, so $E_t$ cannot be Stein.
\end{proof}

\Needspace{10\baselineskip}
\begin{corollary}[= Corollary~\ref{cor:nakano-affine-weak}]
Let $X$ be a compact K\"ahler manifold and let $(E,h)$ be a Hermitian holomorphic vector bundle such that $(E^*,h^*)$ is Nakano semipositive. Then, for every $t\in H^1(X,E)$, the affine bundle $E_t$ is weakly $1$-complete.
\end{corollary}

\begin{proof}
Let $a\in A^{0,1}(X,E)$ be the $\overline\partial$-harmonic representative of $t$, and put $n=\dim_{\C}X$.  Let
\[
 \beta:=\star_h a\in A^{n,n-1}(X,E^*)
\]
be the generalized Hodge dual of $a$, defined by the K\"ahler metric on $X$ and the Hermitian metric $h$ on $E$.  By Hodge--Serre duality, $\beta$ is $\overline\partial$-harmonic.  Denote by $\delta'_{h^*}$ the formal adjoint of the $(1,0)$-part $D'_{h^*}$ of the Chern connection on $E^*$.

Apply the Bochner--Kodaira--Nakano identity to the $E^*$-valued $(n,n-1)$-form $\beta$; see \cite[Chapters~VI--VII]{DemaillyCAG}.  Since $\Delta''\beta=0$,
\[
 0
 =\|D'_{h^*}\beta\|^2
  +\|\delta'_{h^*}\beta\|^2
  +\bigl\langle[\sqrti\Theta_{h^*}(E^*),\Lambda]\beta,\beta\bigr\rangle .
\]
The first term vanishes for bidegree reasons, because $\beta$ has holomorphic degree $n$.  The curvature term is nonnegative by the Nakano semipositivity of $(E^*,h^*)$, hence $\delta'_{h^*}\beta=0$. The Hodge-star formula for the formal adjoint gives, up to an irrelevant sign, $\delta'_{h^*}(\star_h a)=\pm\star_h(D'_h a)$. Therefore $D'_ha=0$.  Since $\overline\partial a=0$ as well, $a$ is $D_h$-parallel.  Finally, Nakano semipositivity of $E^*$ implies Griffiths semipositivity of $E^*$, hence Griffiths seminegativity of $E$.  Theorem~\ref{thm:main-affine} now gives the weak $1$-completeness of $E_t$.
\end{proof}

\begin{example}[Unitary flat affine bundles]\label{ex:unitary-flat-affine}
The affine-bundle criterion admits a simple interpretation when the linear part is unitary flat. Let $(E,h)$ be a unitary flat bundle over a compact K\"ahler manifold $X$, and let $a$ be the harmonic representative of a class in $H^1(X,E)$. The K\"ahler identities and the flatness of $E$ imply $D_ha=0$. On the universal cover, after choosing a flat trivialization, $a$ is a closed $\C^r$-valued $(0,1)$-form. Hence there is an anti-holomorphic map
\[
 F\colon\widetilde X\longrightarrow\C^r
\]
with $dF=a$; its affine equivariance gives the corresponding smooth section of the affine bundle.

If $\omega_{\mathrm{Euc}}$ denotes the standard K\"ahler form on $\C^r$, then
\[
 \omega_{h,a}=-F^*\omega_{\mathrm{Euc}}.
\]
Equivalently, regarding $F$ as a holomorphic map to the conjugate vector space $\overline{\C^r}$ and writing $\omega_{\overline{\C^r}}:=-\omega_{\mathrm{Euc}}$, one has
\[
 \omega_{h,a}=F^*\omega_{\overline{\C^r}}.
\]
Thus Theorem~\ref{thm:main-affine} identifies the Steinness of the affine bundle with the K\"ahler positivity of the cohomology class represented by this pullback form.
\end{example}

We apply Theorem~\ref{thm:main-affine} to canonical extensions. Let $(X,\omega)$ be a K\"ahler manifold. The form $\omega$ can be regarded as an $\Omega_X^1$-valued $(0,1)$-form $a_\omega$, whose Dolbeault class $[a_\omega]=[\omega]\in H^1(X,\Omega_X^1)$ defines the canonical extension $\Omega^1_{[\omega]}\to X$. Let $h_\omega$ be the Hermitian metric on $\Omega_X^1$ induced by $\omega$. Since the K\"ahler form is parallel, $D'_{h_\omega}a_\omega=0$, and with the convention in \eqref{eq:omega-ha} one has $\omega_{h_\omega,a_\omega}=\omega$. Moreover, nonnegative holomorphic bisectional curvature of $\omega$ is equivalent to Griffiths seminegativity of $(\Omega_X^1,h_\omega)$.

We therefore obtain the following.

\begin{corollary}\label{cor:canonical-compact}
Let $(X,\omega)$ be a compact K\"ahler manifold with nonnegative holomorphic bisectional curvature. Then $\Omega^1_{[\omega]}$ is Stein.
\end{corollary}

For the complete case we use the following result of Ni--Tam. Let $(M,g)$ be a complete noncompact Riemannian manifold, fix $o\in M$, let $r(o,\cdot)$ be the distance from $o$, and let $H(x,y,t)$ be the heat kernel. Set
\[
 \zeta(x,t):=\int_M H(x,y,t)\,r(o,y)\,d\mu(y).
\]

\begin{lemma}[Ni--Tam {\normalfont\cite[Proposition~3.1]{NiTam2013}}]\label{lem:nitam}
\leavevmode
\begin{enumerate}
 \item If $\operatorname{Ric}(g)\ge0$, then for every $t>0$ the function $\zeta(\cdot,t)$ is a smooth exhaustion of $M$ and satisfies $|\nabla\zeta|(x,t)\le1$.
 \item If, in addition, $(M,g)$ is K\"ahler and has nonnegative quadratic orthogonal bisectional curvature, then for every $t>0$,
 \[
  \bigl\|\partial\overline\partial\zeta(\cdot,t)\bigr\|_g
  \le \frac{1}{\sqrt{2t}}.
 \]
\end{enumerate}
\end{lemma}

We use only the exhaustion assertion in~(1) and the complex Hessian estimate in~(2). Nonnegative holomorphic bisectional curvature implies both nonnegative Ricci curvature and nonnegative quadratic orthogonal bisectional curvature. Hence, when Lemma~\ref{lem:nitam} is applied to $(X,\omega)$, the function $\psi_t:=\zeta(\cdot,t)$ is exhaustive and, for $t$ sufficiently large, satisfies $\omega+\ddbar\,\psi_t>0$.

\Needspace{12\baselineskip}
\begin{theorem}[= Theorem~\ref{thm:main-complete}]
Let $(X,\omega)$ be a complete K\"ahler manifold with nonnegative holomorphic bisectional curvature. Then the canonical extension $\Omega^1_{[\omega]}$ is Stein.
\end{theorem}

\begin{proof}
If $X$ is compact, the assertion is Corollary~\ref{cor:canonical-compact}. Assume that $X$ is noncompact. Let $\Phi$ be the squared-distance function \eqref{eq:affine-distance} on $\Omega^1_{[\omega]}$. The preceding identities, together with Lemma~\ref{lem:omega-ha}, give
\[
 \ddbar\,\Phi\ge p^*\omega.
\]
Choose $t$ sufficiently large so that $\omega+\ddbar\,\psi_t>0$. Then
\[
\begin{aligned}
 \ddbar\,(\Phi+p^*\psi_t)
 &\ge \sqrti\,h_\omega\bigl(\widetilde D'_{h_\omega}u,\widetilde D'_{h_\omega}u\bigr)\\
 &\quad+p^*(\omega+\ddbar\,\psi_t)>0.
\end{aligned}
\]
The right-hand side is positive definite: the second term is positive on the base, and the first is positive in the vertical directions.
The function $\Phi+p^*\psi_t$ is exhaustive: if a sequence escapes to infinity in the base, the term $p^*\psi_t$ tends to infinity; if its projection remains in a compact subset of $X$ while the sequence escapes in the total space, the fibrewise squared distance $\Phi$ tends to infinity. Hence $\Omega^1_{[\omega]}$ admits a strictly plurisubharmonic exhaustion and is Stein.
\end{proof}

\appendix
\section{Exceptional conormal bundles}\label{app:exceptional-conormal}

In this appendix we prove Proposition~\ref{prop:exceptional-conormal-big}. The argument is presumably well known to experts, but we have not found a reference for the precise statement needed here; we therefore include the proof. The point is to pass from Ancona's ample conormal sheaf for a possibly nonreduced thickening of the exceptional set to a positive symmetric power of the ordinary conormal bundle of the chosen irreducible component.  We use the notion of an ample coherent sheaf recalled before Theorem~\ref{thm:ancona}.  Two elementary permanence properties will be used repeatedly: restriction to a closed analytic subspace and passage to a coherent quotient preserve ampleness.  Both follow immediately from the global-generation definition.

\begin{lemma}\label{lem:appendix-initial-term}
Let $Y$ and $A$ be as in Proposition~\ref{prop:exceptional-conormal-big}, let $\mathcal I=\mathcal I_A$, and put
\[
 E:=N^*_{A/Y}=\mathcal I/\mathcal I^2.
\]
Then there exist an integer $a\ge1$ and a nonzero torsion-free ample coherent sheaf $\mathcal G$ on $A$ together with an injection
\[
 \mathcal G\hookrightarrow\operatorname{Sym}^aE.
\]
In particular, $A$ is Moishezon.
\end{lemma}

\begin{proof}
Let $K$ be the connected component of the exceptional set containing $A$, and let $\mathcal J\subset\mathcal O_Y$ be the ideal supplied by Theorem~\ref{thm:ancona} for $K$.  Thus $V(\mathcal J)_{\mathrm{red}}=K$ and $\mathcal J/\mathcal J^2$ is ample.  We give $A$ and $K$ their reduced structures. Since $V(\mathcal J)_{\mathrm{red}}=K$ and $A\subset K$, every holomorphic function vanishing on $K$ also vanishes on $A$. Hence
\[
 \mathcal J\subset\sqrt{\mathcal J}=\mathcal I_K
 \subset\mathcal I_A=\mathcal I.
\]  Restricting the Ancona conormal sheaf to $A$ gives
\[
 \mathcal F:=(\mathcal J/\mathcal J^2)|_A
 \simeq \mathcal J/(\mathcal J^2+\mathcal I\mathcal J),
\]
and $\mathcal F$ is ample.

On a dense open subset of $A$ meeting no other irreducible component of $K$, the radical of $\mathcal J$ equals $\mathcal I$.  There is therefore a unique integer $a\ge1$ such that, generically along $A$,
\[
 \mathcal J\subset\mathcal I^a,
 \qquad
 \mathcal J\not\subset\mathcal I^{a+1}.
\]
Equivalently, $a$ is the largest generic order of vanishing common to all local sections of $\mathcal J$ along $A$.

We first show, after shrinking $Y$ to a neighbourhood of $A$, that the generic inclusion $\mathcal J\subset\mathcal I^a$ actually holds along all of $A$.  Starting with $\mathcal J\subset\mathcal I$, suppose inductively that $\mathcal J\subset\mathcal I^k$ for some $k<a$.  Since $\mathcal I\mathcal J\subset\mathcal I^{k+1}$, taking the degree-$k$ initial term gives a morphism on $A$
\[
 \sigma_k:\mathcal J/\mathcal I\mathcal J
 \longrightarrow \mathcal I^k/\mathcal I^{k+1}.
\]
Because $A\subset Y$ is a local complete intersection, $E=\mathcal I/\mathcal I^2$ is locally free and the canonical graded morphism gives
\[
 \mathcal I^k/\mathcal I^{k+1}
 \simeq \operatorname{Sym}^kE.
\]
By the definition of $a$, the morphism $\sigma_k$ vanishes on a dense open subset of the irreducible reduced space $A$.  Its image is therefore a torsion subsheaf of the locally free sheaf $\operatorname{Sym}^kE$.  Since a locally free sheaf on an irreducible reduced space is torsion-free, $\sigma_k=0$.  Hence $\mathcal J\subset\mathcal I^{k+1}$, and induction yields $\mathcal J\subset\mathcal I^a$ along $A$.

Now
\[
 \mathcal J^2+\mathcal I\mathcal J\subset\mathcal I^{a+1},
\]
so the degree-$a$ initial term descends to a morphism
\[
 \sigma_a:\mathcal F\longrightarrow
 \mathcal I^a/\mathcal I^{a+1}
 \simeq\operatorname{Sym}^aE.
\]
It is generically nonzero by the choice of $a$.  Let $\mathcal G:=\operatorname{Im}\sigma_a$.  Then $\mathcal G\ne0$; it is torsion-free because it is a subsheaf of the locally free sheaf $\operatorname{Sym}^aE$, and it is ample because it is a quotient of the ample coherent sheaf $\mathcal F$.

Finally, $\operatorname{Supp}\mathcal G=A$.  By Ancona's equivalence between ampleness and cohomological ampleness \cite[Proposition~2.5, p.~92]{Ancona1982}, followed by his Moishezon criterion \cite[Proposition~2.8, p.~93]{Ancona1982}, the support of $\mathcal G$ is Moishezon.  Hence $A$ is Moishezon.
\end{proof}

\begin{lemma}\label{lem:appendix-effective-criterion}
Let $X$ be an irreducible projective variety, let $V$ be a vector bundle on $X$, and let $H$ be an ample line bundle.  If
\[
 H^0\bigl(X,\operatorname{Sym}^mV\otimes H^{-1}\bigr)\ne0
\]
for some $m>0$, then $V$ is big.
\end{lemma}

\begin{proof}
Let
\[
 \pi:\Pp(V^*)\longrightarrow X,
 \qquad
 M:=\mathcal O_{\Pp(V^*)}(1).
\]
By our projectivization convention,
\[
 H^0\bigl(X,\operatorname{Sym}^mV\otimes H^{-1}\bigr)
 \simeq
 H^0\bigl(\Pp(V^*),M^m\otimes\pi^*H^{-1}\bigr).
\]
Choose a nonzero section and let $D$ be its effective zero divisor.  Then
\[
 M^m\simeq \pi^*H\otimes\mathcal O(D).
\]
Since $M$ is relatively ample and $H$ is ample on $X$, there exists an integer $c>0$ such that
\[
 B:=M\otimes\pi^*H^c
\]
is ample on $\Pp(V^*)$.  Consequently
\[
 M^{1+cm}
 \simeq M\otimes(M^m)^c
 \simeq B\otimes\mathcal O(cD).
\]
Thus a positive power of $M$ is the tensor product of an ample line bundle and an effective line bundle.

For completeness, this already implies bigness without invoking a separate form of Kodaira's lemma.  Multiplication by the canonical section of $\mathcal O(qcD)$ gives an injection
\[
 H^0\bigl(\Pp(V^*),B^q\bigr)
 \hookrightarrow
 H^0\bigl(\Pp(V^*),M^{q(1+cm)}\bigr)
\]
for every $q\ge1$.  Since $B$ is ample, its Hilbert polynomial has degree $\dim\Pp(V^*)$ and positive leading coefficient.  The right-hand side therefore has maximal polynomial growth along the subsequence $q(1+cm)$, and $M$ is big.  Hence $V$ is big.
\end{proof}

\begin{proof}[Proof of Proposition~\ref{prop:exceptional-conormal-big}]
Put $E=N^*_{A/Y}$.  By Lemma~\ref{lem:appendix-initial-term}, there are $a\ge1$ and a nonzero torsion-free ample coherent sheaf
\[
 \mathcal G\hookrightarrow\operatorname{Sym}^aE,
\]
and $A$ is Moishezon.  Choose a smooth projective modification
\[
 \mu:\widetilde A\longrightarrow A.
\]
Let $U\subset A$ be a dense open subset such that $\mu$ is an isomorphism over $U$, $\mathcal G|_U$ is locally free, and $\mathcal G|_U\to\operatorname{Sym}^a(E|_U)$ has constant rank.

Let $H$ be an ample line bundle on $\widetilde A$ and put
\[
 \mathcal Q:=\mu_*H^{-1}.
\]
The sheaf $\mathcal Q$ is coherent, and the adjunction morphism
\[
 \mu^*\mathcal Q=\mu^*\mu_*H^{-1}\longrightarrow H^{-1}
\]
is an isomorphism over $\mu^{-1}(U)$.  Since $\mathcal G$ is ample in Ancona's sense, for every sufficiently large $k$ the sheaf
\[
 \operatorname{Sym}^k\mathcal G\otimes\mathcal Q
\]
is generated by global sections.  Pulling those sections back to $\widetilde A$ and composing with the adjunction morphism shows that
\[
 \operatorname{Sym}^k(\mu^*\mathcal G)\otimes H^{-1}
\]
is generated by global sections over $\mu^{-1}(U)$.

The morphism $\mathcal G\to\operatorname{Sym}^aE$ induces, after pullback and multiplication in the symmetric algebra, a morphism
\[
 \operatorname{Sym}^k(\mu^*\mathcal G)
 \longrightarrow \operatorname{Sym}^{ak}(\mu^*E).
\]
This morphism is not identically zero.  Indeed, at a point of $\mu^{-1}(U)$ where the original map is nonzero, choose a vector whose image is a nonzero homogeneous polynomial $p\in\operatorname{Sym}^a(E_x)$; its $k$-th power maps to $p^k\ne0$.  Since the preceding global sections generate the source over $\mu^{-1}(U)$, at least one of them has nonzero image.  Therefore
\[
 H^0\!\left(\widetilde A,
 \operatorname{Sym}^{ak}(\mu^*E)\otimes H^{-1}\right)\ne0.
\]
Lemma~\ref{lem:appendix-effective-criterion} implies that $\mu^*E$ is big.  By the definition of bigness on the compact complex space $A$, this is equivalent to $E=N^*_{A/Y}$ being big.  This proves (1).

Assume now that $\mathcal O_{\Pp(E^*)}(1)$ is semipositive.  Its pullback to
\[
 \Pp((\mu^*E)^*)\longrightarrow\widetilde A
\]
is semipositive and hence nef, so $\mu^*E$ is nef.  By (1), $\mu^*E$ is also big.  Yang's determinant theorem \cite[Proposition~4.8, p.~266]{Yang2017} therefore gives
\[
 \det(\mu^*E)=\mu^*(\det E)
\]
big; Yang's proposition applies to nef and big vector bundles over compact K\"ahler manifolds.  Since $\widetilde A$ is a resolution, this says precisely that $\det E$ is big on $A$.  This proves (2).
\end{proof}

\end{document}